\documentclass[10pt,fleqn,a4paper]{article}

\usepackage{pdfsync}
\usepackage{graphicx}
\usepackage{latexsym}
\usepackage[small,bf]{caption2}
\usepackage{amsfonts}
\usepackage{amssymb}
\usepackage{amsmath}
\usepackage{theorem}
\usepackage{times}
\usepackage{bbm}

\newcommand{\bea}{\begin{eqnarray}}
\newcommand{\eea}{\end{eqnarray}}
\newcommand{\bef}{\begin{figure}}
\newcommand{\enf}{\end{figure}}
\newcommand{\ball}{\begin{array}{ll}}
\newcommand{\bacl}{\begin{array}{cl}}
\newcommand{\bal}{\begin{array}{l}}
\newcommand{\bac}{\begin{array}{c}}
\newcommand{\ea}{\end{array}}

\newcommand{\N}{{\mathbb{N}}}
\newcommand{\R}{{\mathbb{R}}}
\newcommand{\Z}{{\mathbb{Z}}}
\newcommand{\E}{{\mathbb{E}}}
\renewcommand{\a}{{\alpha}}
\renewcommand{\b}{{\beta}}
\renewcommand{\P}{{\mathbb{P}}}
\newcommand{\Lcal}{{\mathcal{L}}}

\def\CQFD{\hfill\hbox{\vrule\vbox to 8pt{\hrule width 8pt\vfill\hrule}\vrule}}

\def\PP{\Bbb{P}}
\def\l{\lambda}

\numberwithin{equation}{section}
\theoremstyle{plain}
\newtheorem{theorem}{Theorem}[section]
\newtheorem{lemma}[theorem]{Lemma}
\newtheorem{corollary}[theorem]{Corollary}
\newtheorem{remark}{Remark}
\newtheorem{theorem2}{Theorem}

\title{Zero-range condensation at criticality}
\author{In\'es Armend\'ariz\thanks{Universidad de San Andr\'{e}s, Vito Dumas 284, B1644BID, Victoria,  Argentina. Email: {\tt iarmendariz@udesa.edu.ar}}, Stefan Grosskinsky\thanks{Mathematics Institute, Zeeman Building, University of Warwick, Coventry CV4 7AL, UK. Email: {\tt S.W.Grosskinsky@warwick.ac.uk}}, Michail Loulakis\thanks{School of Applied Mathematical and Physical Sciences, National Technical University of Athens, 15780 Athens, Greece, and Institute of Applied and Computational Mathematics, FORTH, Heraklion Crete, Greece. Email: {\tt loulakis@tem.uoc.gr}} }

\begin{document}
\sloppy
\maketitle

\begin{abstract}
\noindent
Zero-range processes with decreasing jump rates exhibit a condensation transition, where a positive fraction of all particles condenses on a single lattice site when the total density exceeds a critical value.
We study the onset of condensation, i.e. the behaviour of the maximum occupation number after adding or subtracting a subextensive excess mass of particles at the critical density. We establish a law of large numbers for the excess mass fraction in the maximum, which turns out to jump from zero to a positive value at a critical scale. Our results also include distributional limits for the fluctuations of the maximum, which change from standard extreme value statistics to Gaussian when the density crosses the critical point. Fluctuations in the bulk are also covered, showing that the mass outside the maximum is distributed homogeneously. In summary, we identify the detailed behaviour at the critical scale including sub-leading terms, which provides a full understanding of the crossover from sub- to supercritical behaviour.

\end{abstract}

\vspace{1mm}

\section{Introduction}
The zero-range process is one of the interacting particle systems introduced
in the seminal paper \cite{spitzer70}. The process has unbounded local state
space, i.e. there is no restriction on the number of particles per site, and the
jump rate $g(n)$ at a given site depends only on the number of particles $n$ at that site. This simple zero-range interaction leads to a product structure of the stationary measures \cite{andjel82,spitzer70} and further interest was initially on the existence of the dynamics under general conditions \cite{andjel82} and on establishing hydrodynamic limits. These questions have been successfully addressed in the case of attractive zero-range processes when $g(n)$ is a non-decreasing function, and results are summarized in \cite{kipnislandim}. For such processes with additional space dependence of the rates $g_x$, there is also a number of rigorous results regarding condensation on slow sites \cite{andjeletal00,ferrarisisko07,landim96}.\\

\noindent
More recently, there has been increasing interest in zero-range processes with spatially homogeneous jump rates $g(n)$ decreasing with the number of particles. This results in an effective attraction of the particles and can lead to condensation phenomena. A generic family of models with that property has been introduced in the theoretical physics literature \cite{evans00}, with asymptotic behaviour of the jump rates
\bea\label{r0}
g(n)\simeq 1+\frac{b}{n^\lambda}\quad\mbox{as }n\to\infty\ .
\eea
For $\lambda\in (0,1)$, $b>0$ and for $\lambda =1$, $b>2$ the following phase transition was established using heuristic arguments: If the particle density $\rho$ exceeds a critical value $\rho_c$, the system phase separates into a homogeneous background with density $\rho_c$ and a condensate, a single randomly located lattice site that contains a macroscopic fraction of all the particles. This type of condensation appears in diverse contexts such as traffic jamming, gelation in networks, or wealth condensation in macro-economies, and zero range processes or simple variants have been used as prototype models to explain these phenomena (see \cite{evansetal05} for a review).\\

\noindent
The existence of invariant measures with simple product structure makes the problem mathematically tractable. Jeon, March and Pittel showed in \cite{jeonetal00} that for some cases of zero-range processes the maximum site contains a non-zero fraction of all the particles. Condensation has been established rigorously in \cite{grosskinskyetal03} by proving the equivalence of ensembles in the thermodynamic limit, where the lattice size $L$ and the number of particles $N$ tend to infinity such that $N/L\to\rho$. This implies convergence of finite-dimensional marginals of stationary measures conditioned on a total particle number $N$, to stationary product measures with density $\rho$ in the subcritical case $\rho\leq\rho_c$, and with density $\rho_c$ in the supercritical case $\rho >\rho_c$. In the latter case the condition on the particle number is an atypical event which is most likely realized by a large deviation of the maximum component, and the problem can be described as Gibbs conditioning for measures without exponential moments. It turns out (cf. \cite{armendarizetal08}) that a strong form of the equivalence holds in the supercritical case, which determines the asymptotic distribution of the particles on all $L$ sites. A similar result has been established in \cite{ferrarietal07} on a lattice of fixed size $L$ in the limit $N\to\infty$, and the local equivalence of ensembles result was generalized to processes with several particle species in \cite{grosskinsky08}. More recent rigorous results address metastability for the motion of the condensate \cite{beltranetal10,beltranetal09}.\\

\begin{figure}
\begin{center}
        \includegraphics[width=0.7\textwidth]{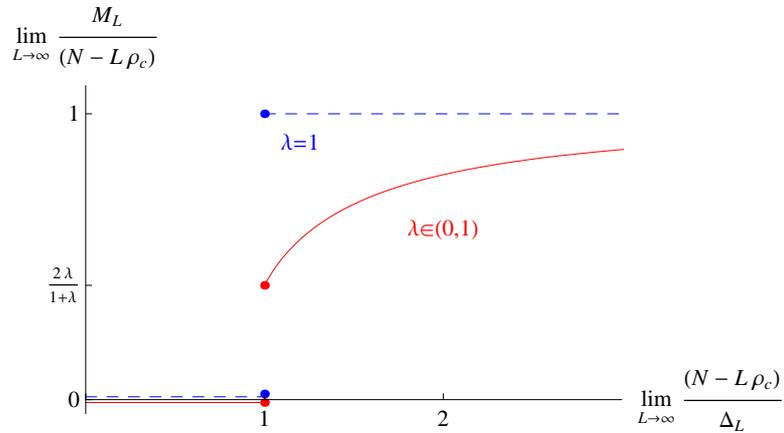}
\end{center}
\caption{   \label{fig:scale}
Illustration of the law of large numbers for the excess mass fraction $\frac{M_L}{N-\rho_c L}$ in the condensate on the critical scale $\Delta_L$ (\ref{cscale0}). For $\lambda =1$ the results are given in Theorem \ref{th1} ((\ref{conv0}) and (\ref{conv})) and for $\lambda\in (0,1)$ in Theorem \ref{th3} ((\ref{conv0l}) and (\ref{convl})). The behaviour at $1$ depends on the sub-leading terms in the excess mass, as detailed in (\ref{gammal1}) and (\ref{subl}).
}
\end{figure}

\noindent
In this paper we study the properties of the condensation transition at the critical density $\rho_c$ for the processes introduced in \cite{evans00} with rates (\ref{r0}), to understand the onset of the condensate formation. We consider the thermodynamic limit with $N/L\to\rho_c$, with the excess mass $N-\rho_c L$ is on a scale $o(L)$. Our results are formulated in Section 2.2 and provide a rather complete picture of the transition from a homogeneous subcritical to condensed supercritical behaviour. It turns out that the condensate forms suddenly on a critical scale $N-\rho_c L\sim \Delta_L$, which is identified in Theorems \ref{th1} and \ref{th3} to be
\begin{equation}\label{cscale0}
\Delta_L =\begin{cases} \sigma\sqrt{(b-3)L\log L}&\mbox{ for}\quad\lambda =1,\ b>3\quad\mbox{and}\\
c_\lambda(\sigma^2 L)^{\frac{1}{1+\lambda}}&\mbox{ for}\quad\lambda \in (0,1),\ b>0\ .\end{cases}
\end{equation}
Our results imply a weak law of large numbers for the ratio $M_L /(N-\rho_c L)$ where $M_L$ is the maximum occupation number, which is illustrated in Figure~\ref{fig:scale}. The ratio exhibits a sudden jump from $0$ to a positive value when the excess mass reaches the critical size $\Delta_L$. At this point both values can occur with positive probability depending on sub-leading orders of the excess mass, which is discussed in detail in Section 2.3. For $\lambda =1$ the full excess mass is concentrated in the maximum right above the critical scale. On the other hand, for $\lambda\in (0,1)$ the excess mass is shared between the condensate and the bulk, and the condensate fraction increases from $2\lambda /(1+\lambda)$ to $1$ only as $(N-\rho_c L)/ \Delta_L \to\infty$. Theorem \ref{bulk} provides results for the bulk fluctuations, which imply that the mass outside the maximum is always distributed homogeneously and the system typically contains at most one condensate site. Theorems \ref{th1} and \ref{th3} also cover the fluctuations of the maximum, which change from standard extreme value statistics to Gaussian. This is complemented by Theorems \ref{th2} and \ref{th4} on downside deviations, which give a detailed description of the crossover to the expected Gumbel distributions in the subcritical regime ($\rho <\rho_c$), where the marginals have exponential tails. In \cite{evansetal08}  the fluctuations of the maximum for $\lambda=1$ were observed by the use of saddle point computations to change from Gumbel ($\rho <\rho_c$), via Fr\'echet ($\rho =\rho_c$), to Gaussian or stable law fluctuations ($\rho >\rho_c$), raising the question on how the transition between these different regimes occurs. Our results around the critical point provide a detailed, rigorous answer to that question, covering also the case $\lambda\in (0,1)$. We use previous results on local limit theorems for moderate deviations of random variables with power-law distribution \cite{doney01} for the case $\lambda =1$, and stretched exponential distribution \cite{nagaev68} for $\lambda\in (0,1)$. In the latter case we can also extend the results for $\rho >\rho_c$ (Corollary \ref{corr}) to parameter values that were not covered by previous results \cite{armendarizetal08}.\\

\noindent
In general, the onset of phase separation and phase coexistence at the critical scale is a classical question of mathematical statistical mechanics. This has been studied for example in the Ising model and related liquid/vapour systems in \cite{biskupetal02,biskupetal03}, where a major point is the shape of critical 'droplets'. Here we treat this question in the case of zero-range condensation, where the main mathematical challenges are related to subexponential scales and a lack of symmetry between the fluid and condensed phase. The condensate turns out to always concentrate on a single lattice site (even at criticality), and contains a positive fraction of the excess mass. In contrast to liquid/vapour systems, this fraction is not 'universal', but depends on the system parameter $\lambda$ (see also discussion in Section 2.3). From a mathematical point of view, the analysis includes interesting connections to extreme value statistics and large deviations for subexponential random variables, which in itself is an area of recent research interest (see \cite{denisovetal08,al10} and references therein). Our results also provide a detailed understanding of finite-size effects and metastability close to the critical point, which are important in applications such as traffic flow and granular clustering (see \cite{chlebounetal10} and references therein).

\section{Definitions and results}

\subsection{The zero-range process and condensation}

We consider the zero-range processes on a finite set $\Lambda_L$ of size $L$. Given a jump rate function $g:\ \N_{0}=\{0,1,2,\ldots\}\mapsto[0,+\infty)$ such that $g(n)=0\Leftrightarrow n=0$ and a set of transition probabilities $p(\cdot,\cdot)$ on $\Lambda_L\times\Lambda_L$, a zero range process is defined as a Markov process on the state space $X_L=\N_{0}^{\Lambda_L}$ of all particle configurations
\bea
\eta =(\eta_x :x\in\Lambda_L )\ ,
\eea
where $\eta_x \in\N_0$ is the local occupation number at site $x$. The dynamics is given by the generator
\bea\label{gen}
\Lcal f(\eta )=\sum_{x,y\in\Lambda_L} g(\eta_x ) p(x,y)\big( f(\eta^{x,y} )-f(\eta )\big)
\eea
using the notation\quad $\eta^{x,y}_z=\left\{\bacl \eta_x -1,\ &\ z=x \,\text{ and }\,\eta_x>0\\ \eta_y +1, &\ z=y\, \text{ and }\,\eta_x>0\\ \eta_z,\ &\mbox{ otherwise.}\ea\right.$ 
\\

\noindent For a technical discussion of the domain of test functions
$f$ of the generator and the corresponding construction of the
semigroup we refer to \cite{andjel82}. The practical meaning of
(\ref{gen}) is that any given site $x$ looses one particle with rate
$g(\eta_x )$ and this particle then jumps to site $y$ with
probability $p(x,y)$. To avoid degeneracies $p(x,y)$ should be
irreducible transition probabilities of a random walk on
$\Lambda_L$. This way, the number of particles is the only conserved
quantity of the process, leading to a family of stationary measures
indexed by the particle density. In the following we are interested
in the situation where these measures are spatially homogeneous.
This is guaranteed by the condition that the harmonic equations \bea
\sum_{x\in\Lambda_L} p(x,y)\lambda_x =\lambda_y\ ,\quad
y\in\Lambda_L, \eea have the constant solution $\lambda_x \equiv 1$,
and by the irreducibility of $p(x,y)$ this implies that every
solution is constant. This is for example the case if $\Lambda_L$ is
a regular periodic lattice and $p(x,y)$ is translation invariant,
such as $\Lambda_L =\Z /L\Z$ and $p(x,y)=\delta_{y,x+1}$ for totally
asymmetric or $p(x,y)=\frac12\delta_{y,x+1} +\frac12\delta_{y,x-1}$ for
symmetric nearest-neighbour hopping.\\

\noindent It is well known (see e.g. \cite{andjel82,spitzer70}) that
under the above conditions the zero-range process has a family of
stationary homogeneous product measures $\nu_\phi$. The occupation
numbers $\eta_x$ are i.i.d. random variables with marginal
distribution
\bea\label{smeas} \nu_\phi \big[\eta_x
=n\big]=\frac{1}{z(\phi )}\, w(n)\,\phi^n \quad\mbox{where}\quad
w(n)=\prod_{k=1}^n \frac{1}{g(k)}\ .
\eea
The parameter $\phi$ of
the stationary measures is called the fugacity, and the measures
exist for all $\phi\geq 0$ such that the normalization (partition
function) is finite, i.e.
\bea
z (\phi ):=\sum_{n=0}^\infty w(n)\,\phi^n <\infty\ .
\eea
The particle density as a function of $\phi$ can be computed as
\bea
R (\phi
):=\E^{\nu_\phi} \big[\eta_x \big]=\phi\ \partial_\phi \log z(\phi
),
\eea
and turns out to be strictly increasing and continuous with $R(0)=0$.\\

\noindent In this paper we consider the family of models introduced
in \cite{evans00}, where the  jump rates have asymptotic behaviour
\bea\label{r0new}
g(n)\simeq 1+\frac{b}{n^\lambda}\quad\mbox{as }n\to\infty\ ,
\eea
with $b>0$ and $\lambda\in (0,1]$. In (\ref{r0new}) and hereafter we  use the notation $a_n\simeq b_n$ as $n\to\infty$, if  $\lim_{n\to\infty}a_n/b_n=1.$ We will also write $a_n\sim b_n$ as $n\to\infty$ if there is a constant $C>1$ such that $C^{-1}\le a_n/b_n\le C$ for sufficiently large $n$. With (\ref{smeas}) this definition of jump rates leads to
stationary weights with asymptotic power law decay
\bea\label{plaw}
w(n)\simeq A_1 n^{-b}\quad\mbox{for }\lambda =1\ ,
\eea
and stretched exponential decay
\bea\label{sexp}
w(n)\simeq A_\lambda \exp\Big(-\frac{b}{1-\lambda}n^{1-\lambda}\Big)\quad\mbox{for
}\lambda \in (0,1)\ ,
\eea
with constant prefactors $A_\lambda$.

In the second case the distributions (\ref{smeas}) are well defined for all $\phi\in [0,1]$ with finite
maximal (\textit{critical}) density
\bea
\rho_c :=R(1)<\infty \label{rhoc} \eea and finite corresponding variance \bea \sigma^2
:=\E^{\nu_1} \big[\eta_x^2 \big]-\rho_c^2 <\infty\ .
\eea
If $\lambda =1$ the corresponding variance is finite if $b>3$, which we
will assume hereafter. The case $2<b\leq 3$ is not covered by our main results, and we discuss it shortly in Section 2.3. In general, (\ref{sexp}) also contains terms of lower order $n^{1-k\lambda}$, $k\geq 2$, in the exponent, which may contribute to the asymptotic behaviour for $\lambda\in (0,1/2]$
depending on the subleading terms in the jump rates (\ref{r0new}). To avoid these complications when  $\lambda\leq 1/2$, we  focus on processes with rates (\ref{r0new}) for which (\ref{sexp}) holds. The simplest way to meet this condition is to choose $g(n)=w(n-1)/w(n)$, $n\ge 1$, with $w(n)$ as in the right hand of (\ref{sexp}) with $A_\lambda=1$.\\

\noindent
It has been shown in
\cite{evans00,grosskinskyetal03} that when the critical density is
finite the system exhibits a condensation transition that can be
quantified as follows. Since the number of particles is conserved by
the microscopic dynamics for each $N\in\N$, the subspaces
\bea
X_{L,N} =\big\{\eta\in X_L :S_L (\eta )=N\big\}\quad\mbox{where}\quad S_L (\eta )=\sum_{x\in\Lambda_L} \eta_x
\eea
are invariant. The zero range process is irreducible
over each of these subspaces and the unique invariant measure
supported on $X_{L,N}$ is given by
\bea
\mu_{L,N} =\nu_\phi \big[\cdot\ |S_L =N\big]\ .
\eea
It is not hard to see  that the
measures $\mu_{L,N}$ are independent of $\phi$ on the right-hand
side. A question of interest is the convergence of the measures
$\mu_{L,N}$ in the thermodynamic limit $L,N\to\infty$, $N/L\to\rho$.
This is answered by the equivalence of ensembles principle, which
states that in the limit the measures $\mu_{L,N}$ locally behave
like a product measure $\nu_\phi$ for a suitable $\phi$. Note that
when $\rho\le\rho_c$ there exists a unique $\phi =\phi(\rho)$ such that
$\rho=R\big(\phi\big)$, whereas if $\rho>\rho_c$ no such
$\phi$ exists. The equivalence of ensembles precisely states that if
$f$ is a cylinder function, i.e. a function that only depends on the
configuration $\eta$ on a finite number of sites, then
\bea\label{equi}
\mu_{L,N} \big[f\big]\to \nu_\phi \big[f\big]\ ,
\eea
provided that (see \cite{grosskinsky08} and Appendix 2.1 in
\cite{kipnislandim})
\bea
R(\phi )=\rho\quad\mbox{and}\quad f\in
L^2 (\nu_\phi )\quad &\mbox{for}
&\quad\rho <\rho_c \quad\mbox{and}\nonumber\\
\phi =\phi_c =1\quad\mbox{and}\quad f\mbox{ bounded}\quad
&\mbox{for} &\quad\rho \geq\rho_c \ .
\eea
The behaviour described
above is accompanied by the emergence of a condensate, a site which
contains $O(L)$ particles. If $\rho<\rho_c$ one can easily check
that the limiting measures $\nu_{\phi(\rho)}$ have finite
exponential moments and the size of the maximum component $M_L (\eta
)=\max_{x\in\Lambda_L} \eta_x$ is typically $O(\log L)$. If on the
other hand $\rho >\rho_c$ it has been shown in \cite{jeonetal00} for
the power law case that
\bea\label{lln}
\frac1L M_L
\stackrel{\mu_{L,N}}{\longrightarrow}\rho -\rho_c\ .
\eea
The notation in (\ref{lln}), which we also use in the following, denotes
convergence in probability w.r.t the conditional laws $\mu_{L,N}$,
i.e.
\bea
\mu_{L,N} \Big[\ \Big|\frac1L M_L -(\rho -\rho_c )\Big|
>\epsilon\Big]\to 0\quad\mbox{for all }\epsilon >0\ .
\eea
Equation (\ref{equi}) has been generalized in \cite{armendarizetal08} for
$\rho >\rho_c$ to test functions $f$ depending on all sites but the
maximally occupied one, and equation (\ref{lln}) is proved for
stretched exponentials of the form (\ref{sexp}) with $\lambda >1/2$,
as well. An immediate corollary is that the size of the second
largest component is typically $o(L)$, which implies that the
condensate typically covers only a single randomly located 
site.

\subsection{Main results}

In the following we study the distribution of the excess mass in the
system at the critical point to fully understand the emergence of the condensate when the density increases from sub- to supercritical values. We consider the thermodynamic limit
$N/L\to\rho_c$ where the excess mass is on a sub-extensive scale $|N-\rho_c
L|=o(L)$. Our first theorem on the power-law case (\ref{plaw})
relies on a result of Doney (Theorem 2 in \cite{doney01}) for the
estimation of $\nu_{\phi_c}\big[S_L=N\big]$. Precisely, for
$z:=(N-\rho_c L )/\sqrt{L}\to\infty$, we get
\begin{equation}
\nu_{\phi_c} \big[S_L =N\big]=\frac{1}{\sqrt{2\pi\sigma^2 L}}\, e^{-\frac{z^2}{2\sigma^2}}\big( 1+o(1)\big) +L\ \nu_{\phi_c}\Big[\eta_x=\big[ z\sqrt{L}\big]\Big]\big( 1+o(1)\big)
\label{split}
\end{equation}
as $L\to\infty$. 
It turns out that when $N-\rho_cL$ is close to the typical scale $\sqrt{L}$ (case a) in the Theorem below) the first term of the sum dominates the right hand side of (\ref{split}) and the excess mass is distributed homogeneously among the sites. On the other hand, when $N-\rho_cL$ is large enough (case (b)) it is the second term that dominates the right hand side of (\ref{split}) and this implies the existence of a condensate that carries essentially all the excess mass. Finally, there is an intermediate scale (case (c)) where the two terms are of the same order and both scenarios can occur with positive probability.

\begin{theorem}
{\bf (Upside moderate deviations, power law case)}

\noindent Let $\lambda =1$ and $b>3$, so that $\sigma^2<\infty$.
Assume that $N\ge\rho_c L$
and define $\gamma_L \in\R$ by
\begin{equation}
N=\rho_c L+\sigma\sqrt{(b-3)L\log
L}\,\Big(1+\frac{b}{2(b-3)}\frac{\log\log L}{\log
L}+\frac{\gamma_L}{\log L}\Big). \label{gammal1}
\end{equation}
a) If $\gamma_L\to-\infty $ the distribution under $\mu_{L,N}$ of the maximum $M_L$ is asymptotically equivalent to its distribution under $\nu_{\phi_c}$. Precisely, for all $x>0$ we have
\bea\label{mda}
\lim_{L\to\infty}\mu_{L,N}\left[ \frac{M_L}{L^{\frac{1}{b-1}}}\le x\right]=\lim_{L\to\infty}\nu_{\phi_c}\left[ \frac{M_L}{L^{\frac{1}{b-1}}}\le x\right]=\exp\left\{-\frac{A_1 x^{1-b}}{b-1}\right\}. 
\eea

\begin{equation}\label{conv0}
\text{In particular,\ if }\ N-\rho_cL\gg L^{\frac{1}{b-1}} \ \text{ then }\qquad
\frac{M_L}{N-\rho_c L} \stackrel{\mu_{L,N}}{\longrightarrow}0.
\end{equation}
b) If $\gamma_L\to+\infty $ the normalized fluctuations of the maximum around the excess mass under $\mu_{L,N}$ converge in distribution to a normal r.v.,
\bea\label{fluct}
\frac{M_L -(N-\rho_c L )}{\sqrt{L\sigma^2}} \stackrel{\mu_{L,N}}{\Longrightarrow} {\cal N} (0,1)\ .
\eea
\bea\label{conv}
\text{In particular, }\qquad
\frac{M_L}{N-\rho_c L} \stackrel{\mu_{L,N}}{\longrightarrow}1.
\eea
c) If $\gamma_L \to\gamma\in\R$ we have convergence in distribution
to a Bernoulli random variable, \bea\label{bern} \qquad\ \qquad\
\qquad\ \frac{M_L}{N-\rho_c L} \stackrel{\mu_{L,N}}{\Longrightarrow} 
Be(p_\gamma )\ , \eea where $p_\gamma \in (0,1)$ is such that
$p_\gamma \to 0\ (1)$, as $\gamma\to -\infty\ (+\infty)$. An
explicit expression for $p_\gamma$ is given in (\ref{lgdef}) and (\ref{bernoulimit})
in the proofs section. \label{th1}
\end{theorem}
\noindent
The next result connects the fluctuations of the maximum to the extreme value statistics expected in the subcritical regime.
 
\begin{theorem}
{\bf (Downside moderate deviations, power law case)}

\noindent Let $\lambda =1$ and $b>3$ and define $\omega_L \geq 0$ by
\bea
N=\rho_c L-\omega_L \sigma^2 L^{\frac{b-2}{b-1}}\ .
\eea
\noindent a) If 
$\omega_L \to 0$ the distribution under
$\mu_{L,N}$ of the maximum $M_L$ is asymptotically equivalent to its
distribution under $\nu_{\phi_c}$. Precisely, for all $x>0$ we get
\bea\label{1dpl} \lim_{L\to \infty}\mu_{L,N}\left[
\frac{M_L}{L^{\frac{1}{b-1}}}\le x\right]=\exp\left\{-\frac{A_1 x^{1-b}}{b-1}\right\}. \eea
b) If 
$\omega_L \to \omega >0$ then there exists a positive
constant $v$ such that for all $x>0$ \bea\label{2dpl} \lim_{L\to
\infty} \mu_{L,N}\left[\frac{M_L}{L^{\frac{1}{b-1}}}\le x
\right]=\exp\left\{-A_1 \int_x^{\infty}e^{-\omega
t}\frac{dt}{t^b}\right\}.\eea
c) If 
$\omega_L \to \infty$ then there exist sequences $B_L\to \infty$ and
$s_L= \frac{\rho_cL-N}{\sigma^2 L}(1+o(1))$ with $B_Ls_L\to \infty$, such that for all
$x\in\mathbb{R}$ \bea\label{3dpl} \lim_{L\to \infty} \mu_{L,N}\left[\frac{M_L-B_L}{1/s_L}\le x\right]=\exp\{-e^{-x}\}.\eea
\label{th2}
\end{theorem}

\noindent We return to a more detailed discussion of these results in Section 2.3 after stating the
results for the stretched exponential tail ($\lambda<1$). For this case the
counterpart of estimate (\ref{split}) was obtained by A.V. Nagaev in
\cite{nagaev68}, where the size of the maximum is also discussed and which is summarized in the appendix. In
fact, a careful reading reveals that equation (\ref{fluct2}) below
is already contained there.


\begin{theorem} {\bf (Upside moderate deviations, stretched exponential case)}

\noindent Let $\lambda \in (0,1)$ and $c_\lambda=(1+\lambda)(2\lambda)^{-\frac{\lambda}{1+\lambda}}\left(\frac{b}{1-\lambda}\right)^{\frac{1}{1+\lambda}}$.  Assume that $N\ge\rho_cL$
and define $t_L \geq 0$ by
\bea\label{gammal2}
N=\rho_c L+t_L(\sigma^2 L)^{\frac{1}{1+\lambda}}.
\eea
a) If $\lim\sup t_L<c_\lambda $ the distribution under $\mu_{L,N}$ of the maximum
$M_L$ is asymptotically equivalent to its distribution under $\nu_{\phi_c}$. Precisely,
there exist sequences $y_L,b_L$ such that
\bea
y_L\simeq\Big(\frac{1-\lambda}{b}\log L\Big)^{\frac{1}{1-\lambda}},\qquad b_L\simeq\frac{y_L^\lambda}{b} \qquad\text{as } L\to\infty,
\label{extremes}
\eea
\noindent
and for all $x\in\R$ we have
\bea
\lim_{L\to\infty}\mu_{L,N}\left[ \frac{M_L-y_L}{b_L}\le x\right]=\lim_{L\to\infty}\nu_{\phi_c}\left[ \frac{M_L-y_L}{b_L}\le x\right]=e^{-e^{-x}}.
\label{mdal}
\eea
\begin{equation}
\text{In particular, \ if }\ N-\rho_c L\gg (\log L)^{\frac{1}{1-\lambda}} \ \text{ then}\qquad
\frac{M_L}{N-\rho_c L} \stackrel{\mu_{L,N}}{\longrightarrow}0.
\label{conv0l}
\end{equation}
b) If $t_L\to t$ with $c_\lambda<t\le +\infty$, there exists a
sequence $a_L$ and a function  $a(t)$, $a_L\to a(t)$, such that
\bea\label{fluct2}
\frac{M_L-(N-\rho_cL)a_L}{\sqrt{L}}\stackrel{\mu_{L,N}}{\Longrightarrow}
{\cal N}\left(0,
\frac{\sigma^2}{1-\frac{\lambda\big(1-a(t)\big)}{a(t)}}\right).
\eea
\bea\label{convl} \text{In particular, }\qquad \frac{M_L}{N-\rho_c
L} \stackrel{\mu_{L,N}}{\longrightarrow}a(t).
\eea
The sequence $a_L$ is implicitly defined by (\ref{alef}) in the Appendix (with $a_L =1-\alpha$ and $\gamma =b/(1-\lambda)$), and when $N-\rho_c L\gg L^{\frac{1}{2\lambda}}$ we may take $a_L=1$. The
limit $a(t)$ in the preceding equation is an increasing function of
$t$ with 
\bea\label{nonuni}
\lim_{t\downarrow
c_\lambda}a(t)=\frac{2\lambda}{1+\lambda}\qquad\text{and}\qquad\lim_{t\uparrow\infty}a(t)=a(+\infty)=1.
\eea
c) If $t_L \to c_\lambda$, assume that $\lambda >1/2$ and suppose
\begin{equation}
N=\rho_c L+ c_\lambda(\sigma^2 L)^{\frac{1}{1+\lambda}}-\frac{1+\lambda}{2\lambda c_\lambda}(\sigma^2 L)^{\frac{\lambda}{1+\lambda}}\big(\tfrac{3}{2}\log L +\gamma_L \big)
\label{subl}
\end{equation}
with $\gamma_L\to \gamma\in\R$. Then we have convergence to a
Bernoulli random variable, 
\bea
\frac{M_L}{N-\rho_c L}\stackrel{\mu_{L,N}}{\Longrightarrow} 
\frac{2\lambda}{1+\lambda} Be
(p_\gamma )\ , \label{bern2} \eea
where $p_\gamma \in (0,1)$ is such that $p_\gamma \to 0\ (1)$ for $\gamma\to -\infty\ (+\infty)$. An explicit expression for $p_\gamma$ is given in (\ref{pstretched}).
\label{th3}
\end{theorem}
In c) analogous statements also hold for the case $\lambda\leq 1/2$, which can be derived from the results in \cite{nagaev68} summarized in the appendix. However, the order of the sub-leading scale depends on the first few Cram\'er coefficients of the distribution, and results cannot be formulated in an explicit form as above.

\begin{theorem} {\bf (Downside moderate deviations, stretched exponential case)}

\noindent Let $\lambda<1$ and define $\omega_L\ge 0$ by \[N=\rho_cL-\omega_LL\ (\log L)^{-\frac{1}{1-\lambda}}.\]\\
If $\omega_L\to c\in[0,+\infty]$
there exist sequences $\gamma_L$ and $\zeta_L$, both increasing to $\infty$ with $L$, such that
\bea
\lim_{L\to \infty}\mu_{L,N}\left[ M_L\le \gamma_L+x\,\zeta_L\right]=e^{-e^{-x}},\quad x\in \R\,.
\label{sexponential-limit}
\eea
If $\omega_L\to 0$ the distribution under $\mu_{L,N}$ of the maximum $M_L$ is asymptotically equivalent to its distribution under $\nu_{\phi_c}$. Precisely, if $y_L$ and $b_L$ are the sequences introduced in Theorem \ref{th3}.a), we can take $\gamma_L=y_L$ and $\zeta_L=b_L$ to get
\[
\lim_{L\to \infty}\mu_{L,N}\left[ M_L\le y_L+x\, b_L\right]=e^{-e^{-x}},\quad x\in \R\,.
\]

\label{th4}
\end{theorem}

\noindent Our final result focuses on the fluctuations of the bulk outside the maximum.

\begin{theorem} ({\bf Fluctuations of the bulk})\\
Assume $\lambda\in(0,1)$, or $\lambda=1$ and $b>3$ so that $\sigma^2<+\infty$. \\
a) In the subcritical regime, that is if $\frac{M_L}{\Delta_L}\stackrel{\mu_{L,N}}{\longrightarrow} 0$
, the distribution under $\mu_{L,N}$ of the bulk fluctuation process
converges in the Skorokhod space to a standard Brownian bridge conditioned to return to the origin at time 1, i.e.
\[
X_s^L=\frac{1}{\sigma\sqrt L}\sum_{x=1}^{[sL]} \left(\eta_x-\frac{N}{L}\right)\stackrel{\mu_{L,N}}{\Longrightarrow} BB_s\ .
\]
b) In the supercritical regime, that is if $N-\rho_c L>\Delta_L$ and $\frac{M_L}{N-\rho_cL}\stackrel{\mu_{L,N}}{\longrightarrow} \kappa$, $\kappa$ a positive constant, the distribution under $\mu_{L,N}$ of the bulk fluctuation process, converges in the Skorokhod space to a standard Brownian bridge plus an independent, random drift term. Precisely, if $\tilde{\eta}_x=\eta_x\mathbbm{1}\{\eta_x\le L^{1/4}\}$ then
\[
Y_s^L=\frac{1}{\sigma\sqrt L}\sum_{x=1}^{[sL]} \left(\tilde{\eta}_x-\frac{N-a_L(N-\rho_cL)}{L}\right)\stackrel{\mu_{L,N}}{\Longrightarrow} 
BB_s+s\,\Phi\,,
\]
where $\Phi\sim {\cal N}\Big(0,\,1/\big(1-\frac{\lambda(1-a(t))}{a(t)}\big)\Big)$\quad is an independent random variable.\\
When $\lambda=1$, or when $\lambda\in(0,1)$ and $N-\rho_cL\gg L^{\frac{1}{2\lambda}}$ we may take $a_L=1$. Otherwise, $a_L$ is defined by (\ref{alef}) in the Appendix (with $a_L =1-\alpha$ and $\gamma =b/(1-\lambda )$).
\label{bulk}
\end{theorem}

\noindent
The supercritical case (assertion b) above) takes a particularly simple form when $a(t)=1$: then $\Phi\sim{\cal N}(0,1)$ is a Gaussian variable independent of the Brownian bridge component, and hence $BB_s+s\,\Phi$ is a standard Brownian motion $B_s$. This is the case for the supercritical power law, or for the stretched exponential law when $\frac{N-\rho_cL}{\Delta_L}\to +\infty$.

\subsection{Discussion of the main results}

\begin{figure}
\begin{center}
        \includegraphics[width=0.7\textwidth]{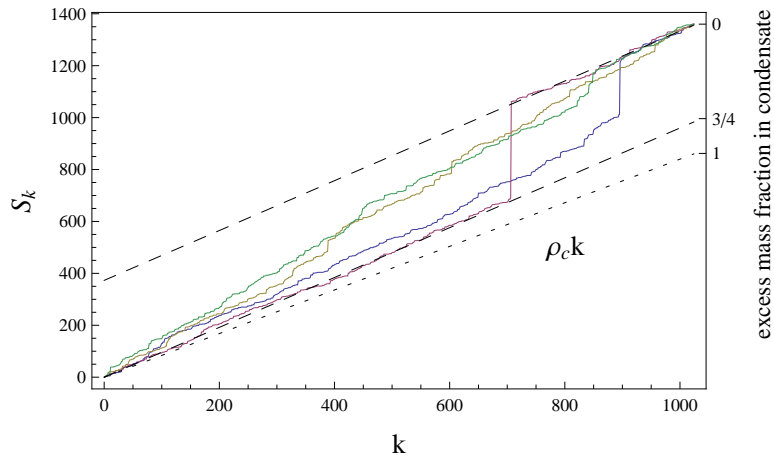}
\end{center}
\caption{   \label{fig:profile}
Results of Monte Carlo simulations of the zero-range process with rates (\ref{r0}) with $\lambda =0.6$ and $b=2$, on a one-dimensional lattice with $L=1024$ sites and periodic boundary conditions. For these parameter values $\rho_c =0.842$, $\sigma^2 =2.55$ and $c_\lambda =4.09$, and we choose $N=1360$ particles, which is very close to the leading order prediction of $1356$ for the critical scale according to (\ref{subl}) with $\gamma =0$. We plot four realizations for the accumulated profile $S_k=\sum_{x=1}^k \eta_x$ against $k$, and see that both, fluid and condensed realizations occur. In the condensed case, the mass is shared between the condensate and the bulk according to the prediction (\ref{nonuni}) with $\frac{2\lambda}{1+\lambda}=3/4$, as is indicated by the dashed lines.
}
\end{figure}

As is already summarized in the introduction, Theorems \ref{th1} and \ref{th3} imply a weak law of large numbers for the excess mass fraction in the condensate $M_L /(N-\rho_c L)$. The critical scale $\Delta_L$ for the excess mass, above which a positive fraction of it concentrates on the maximum and forms a condensate according to (\ref{conv0}), (\ref{conv}) and (\ref{conv0l}), (\ref{convl}), is summarized in (\ref{cscale0}). It is of order $\sqrt{L\log L}$ for the power law case given precisely in (\ref{gammal1}), and the lighter tails in the stretched exponential case lead to a higher scale of order $L^{\frac{1}{1+\lambda}}$ given precisely in (\ref{gammal2}). At the critical scale the excess mass fraction can take both values with positive probability (cf. (\ref{bern}) and (\ref{bern2})), depending on sub-leading orders as detailed in (\ref{gammal1}) and (\ref{subl}). In the power law case, the condensate always contains the full excess mass (\ref{conv}) as soon as it exists. On the other hand, for stretched exponential tails the excess mass is shared between the condensate and the bulk according to (\ref{convl}) as long as $N-\rho_c L\sim \Delta_L$, and the fraction $a(t)$ of the condensate gradually increases to $1$. This behaviour is illustrated in Figure \ref{fig:scale} in the introduction. The results on the bulk fluctuations in Theorem \ref{bulk} imply that below criticality the excess mass is distributed homogeneously in the system, and that the same holds above criticality in the bulk, which completes the above picture. These results are illustrated in Figure \ref{fig:profile}, where we show sample profiles for a zero-range process which show exactly the predicted behaviour already at a rather moderate system size of $L=1024$.\\

\noindent
The discontinuous formation of the condensate on the critical scale implies that it forms 'spontaneously' out of particles taken from the bulk of the system: When crossing the critical scale by adding more mass to the system, the number of particles joining the maximum is indeed of higher order than the number of particles that have to be added to the system in order to form the condensate. A similar phenomenon has been reported for the Ising model and related liquid/vapour systems in \cite{biskupetal02,biskupetal03}. In contrast to these results, the condensed excess mass fraction at criticality is not 'universal', but depends on the system parameter $\lambda$ according to (\ref{nonuni}). This might seem surprising at first sight, but the rates of the form (\ref{r0}) introduce an effective long-range interaction when the zero-range process is mapped to an exclusion model with finite local state space (see e.g. \cite{evans00,evansetal05}).\\



\noindent
In addition to a law of large numbers our results also include limit theorems for the fluctuations of the maximum, which are Gaussian above the critical scale (cf. (\ref{fluct}) and (\ref{fluct2})), and given by the extreme value statistics below criticality. As long as $\lim_{L\to\infty} (N-\rho_c L)/ \Delta_L <1$, the excess mass does not affect the behaviour of the maximum. According to statements a) of Theorems \ref{th1} and \ref{th3}, $M_L$ scales as the maximum of i.i.d.\ random variables, which is proportional to $L^{\frac{1}{b-1}}$ and $(\log L)^{\frac{1}{1-\lambda}}$, respectively, with limiting Fr\'echet distribution for power law tails (\ref{mda}) (cf. \cite{evansetal08}) and Gumbel distribution for stretched exponential tails (\ref{mdal}). Theorems \ref{th2} and \ref{th4} describe the crossover to the expected Gumbel distributions in the subcritical regime, where the marginals have exponential tails. In the power law case, the change from Fr\'echet to Gumbel occurs at the critical scale $\rho_c L-N\sim L^{(b-2)/(b-1)}$ according to (\ref{2dpl}). In \cite{evansetal08} the behaviour of the maximum was predicted for $N=\rho L$ with $\rho$ smaller, equal, or bigger than $\rho_c$ for the power law case $\lambda =1$. Our results provide a rigorous confirmation including the stretched exponential case $\lambda\in (0,1)$, together with a full understanding of the crossover from subcritical extremal statistics to Gaussian fluctuations in the supercritical regime.\\

\noindent
We point out that at criticality the correlations introduced by conditioning on the total number of particles shift from being entirely absorbed by the bulk to being entirely absorbed by the maximum. Indeed, when $N-\rho_c L\ll \Delta_L$ we know from Theorems \ref{th1} a) and \ref{th3} a) that the maximum behaves as the maximum of i.i.d. random variables with distribution $\nu_{\phi_c}$. On the other hand, if $N-\rho_c L\gg \Delta_L$, the bulk asymptotically behaves as i.i.d. random variables with distribution $\nu_{\phi_c}$ following from Theorems 1a and 1b in \cite{armendarizetal08}, and the discussion after Theorem \ref{bulk}.\\

\noindent
In the stretched exponential case, there is another interesting point regarding the centering of the bulk variables in the central limit theorem: 
When the excess mass exceeds $\Delta_L$
the typical excess mass in the bulk is 
\[
(1-a_L)(N-\rho_c L)\sim\frac{\sigma^2 bL}{(N-\rho_cL)^\lambda}\ ,
\]
as follows from the implicit definition (\ref{alef}) of $a_L$. This is of order at least $\sqrt{L}$ unless $N-\rho_c L\gg L^{1/(2\lambda )}$, hence the special centering required in Theorem \ref{bulk} b). In this case
 the equivalence of ensembles cannot be extended to the strong form of  \cite{armendarizetal08} (Theorem 2b).
Note that for $\lambda \leq 1/2$ this affects even supercritical densities, i.e. $N/L\to\rho>\rho_c$. This is why previous results did not cover this case, which is summarized in the following simple Corollary of Theorem \ref{th3}, and completes the condensation picture for supercritical densities.

\begin{corollary}
If $\lambda\in (0,1/2]$ and $N/L\to\rho >\rho_c\ $ we have
\[
\frac{1}{L}\, M_L \stackrel{\mu_{L,N}}{\longrightarrow} \rho-\rho_c\,,
\]
and the fluctuations around this limit are given by (\ref{fluct2}), with $(1-a_L)=O(L^{-\lambda})$.
\label{corr}
\end{corollary}

\noindent
A necessary condition for our results is the existence of finite second moments, and the case $\lambda=1$ and $2<b\le 3$ is not covered by this article. The reason we cannot provide results analogous to Theorems \ref{th1} and \ref{th2} is the lack of a precise estimate for the probability of a moderate deviation of the sum in that case, similar to the result 
(\ref{split}) by Doney \cite{doney01} for square integrable power-law tails. Nevertheless, when the excess mass is such that
\[
\P\big[S_L=N \big]=L\,\P\Big[\eta_1=[N-\rho_c L]\Big] \big(1+o(1)\big),
\]
we can still apply Theorem 1 in \cite{armendarizetal08} to obtain a stable limit theorem for the fluctuations of the maximum around $N-\rho_cL$. For instance if $2<b<3$, the preceding relation is true provided $N-\rho_cL{\gg}L^{\frac{1}{b-1}}$, and under this condition we get that
\[
\frac{M_L}{N-\rho_cL}\stackrel{\mu_{L,N}}{\longrightarrow}1\quad\mbox{and}\quad
\frac{M_L-(N-\rho_cL)}{L^{\frac{1}{b-1}}}\stackrel{\mu_{L,N}}{\Longrightarrow}G_{b-1}\ ,
\]
where $G_{b-1}$ is a completely asymmetric stable law with index $b-1$.

\section{Proofs}

Since the product measures $\nu_\phi$ and the conditional distributions $\mu_{L,N}$ are exchangeable and independent of the jump probabilities $p(x,y)$, the spatial structure of zero-range configurations is irrelevant for our results. In the following we will therefore consider $\eta_1 ,\eta_2 ,\ldots$ to be i.i.d. integer valued random variables defined in a probability space $(\Omega,{\cal F},\P)$, where ${\cal F}$  is the $\sigma$-field generated by $\{\eta_i\}_{i\in\N}$ and $\P=\nu_{\phi_c}$. We further define
\bea p_k:=\P \big[\eta_i =k\big]. \eea
Note that $p_n$ is directly proportional to the stationary weights
$w(n)$ in (\ref{smeas}). Recall the notation
\bea\label{notation}
\rho_c =\E\big[\eta_i \big]\ ,\quad\sigma^2 =\E \big[\eta_i^2
\big]-\rho_c^2 \ ,\quad S_L =\sum_{i=1}^L \eta_i \ ,\quad M_L
=\max_{i=1,..,L} \eta_i \ , \eea and that the conditional laws are
given by $\mu_{L,N} =\P \big[\cdot\ |\ S_L =N\big]$. We will denote by $x^+=\max\{x,0\}$ the positive part of a real number $x$, and by $x^-=(-x)^+$ its negative part.\\

\subsection{Preliminaries}

\noindent Our proofs mainly involve explicit estimates and standard
large deviations methods. One such technique consists in introducing
a change of measure that renders the rare event typical. Precisely,
given $\alpha>0$ and $s\in \R$, define a new measure
$\P_{\alpha}(s)$ on the $\sigma$--field ${\cal
F}_L:=\sigma(\eta_1,\dots,\eta_L)$
 by
\begin{equation} \frac{d\P_{\alpha}(s)}{d\P}\Big|_{{\cal F}_L}=
\frac{1}{Z^L_{\alpha}(s)}\mathbbm{1}_{\{M_L\le \alpha\}}e^{s
S_L}\,,
\label{tilted}
\end{equation} 
where the normalization above is given by
\[ Z_{\alpha}(s)=\sum_{k\le \alpha} e^{sk}p_k\,.  \]
Note that under $\P_{\alpha}(s)$ the random variables
$\{\eta_k\}_{k\in\N}$ are i.i.d., bounded above by $\alpha$, their
mean value is given by
\[
\rho_{\alpha}(s)=\frac{Z'_{\alpha}(s)}{Z_{\alpha}(s)}=\frac{1}{Z_{\alpha}(s)}
\sum_{k\le\alpha}ke^{sk}p_k,
\]
and their variance is given by
$\sigma_{\alpha}^2(s)=\rho'_{\alpha}(s)$. It is not hard to verify that 
\[
\lim_{s\to-\infty}\rho_\alpha(s)=\inf\{k\ge 0: p_k>0\} \qquad\text{and}\qquad \lim_{s\to\infty}\rho_\alpha(s)=\sup\{k\le\alpha: p_k>0\}.
\]
Since $\rho_\alpha(\cdot)$ is a continuous increasing function, it follows that if $N/L$ is sufficiently close to the mean $\rho_c$ of the
distribution and $\alpha$ is sufficiently large, there exists an $s_*=s_*(L,N,\alpha)$ such that 
\begin{equation}
\rho_\alpha(s_*)=\frac{N}{L}.
\label{sstar}
\end{equation}
The following lemma can be applied to compute the exact asymptotics of the conditional maximum when the average is set to be a small perturbation of the mean, using an a priori estimate as input.
\begin{lemma}
Take $N=N(L)$ such that $\frac{N}{L}\to\rho_c$ and suppose the following conditions are satisfied:\\

\noindent
i) There exists a sequence $\alpha_{L,N}\le\infty$ such that  \[\lim_{L\to\infty}\mu_{L,N} \big[M_L\le\alpha_{L,N}\big]=1\qquad\text{and}\qquad  \lim_{L\to\infty}\ \left(\frac{N}{L}-\rho_c\right)^+\alpha_{L,N}=0.\]

\noindent
ii) For each $x\in\mathbb{R}$, there exists a sequence $\beta_L(x)\leq\alpha_{L,N}$ with $\beta_L(x)\to+\infty$, such that
\[ L\sum_{\beta_L(x)< k\le \alpha_{L,N}} e^{s_*k} p_k \longrightarrow \Phi(x)\in [0,\infty )\quad\text{as}\quad L\to\infty\, ,\]
where $s_*=s_*(L)$ is defined as in (\ref{sstar}) with $\alpha=\alpha_{L,N}$\,. 
Then
\[
\mu_{L,N}\left[M_L\le \b_L(x)\right] \longrightarrow e^{-\Phi(x)}\qquad
\text{as}\quad L\to \infty\, .
\]
\label{lemma0}
\end{lemma}
%
%
%
%
{\em Proof:} We begin by showing that under the conditions of the Lemma $s_*^+\alpha_{L,N}\to 0.$
For ease of notation we may write $\alpha$ and $\beta$ as shorthands for $\alpha_{L,N}$ and $\beta_L(x)$ respectively. Using the elementary inequality $(x-y)(e^x-e^y)\ge(x-y)^2$, valid for all $x,y\ge 0$ we have for any $s>0$
\begin{align}
&Z_{\alpha}(s)\left(\rho_{\alpha}(s)-\frac{N}{L}\right) =\sum_{k\le\alpha}\big(k-\frac{N}{L}\big)e^{sk}p_k\nonumber\\
&\qquad\ge e^{sN/L}\sum_{k\le\alpha} \big(k-\frac{N}{L}\big)p_k+s\sum_{k\le \alpha} \big(k-\frac{N}{L}\big)^2p_k\nonumber\\
&\qquad\qquad\ge e^{sN/L}\left(\rho_c-\frac{N}{L}-\sum_{k>\alpha}\big(k-\frac{N}{L}\big)p_k\right)+s\sigma^2-s\sum_{k>\alpha}\big(k-\frac{N}{L}\big)^2p_k.
\label{s0}
\end{align}
If we set
\[
s_0\sigma^2=2\sum_{k>\alpha}\big(k-\frac{N}{L}\big)p_k+2\big(\frac{N}{L}-\rho_c\big)^+
\]
it follows from (\ref{s0}) that $\rho_\alpha(s_0)>N/L$ for sufficiently large $L$, and since $\rho_\alpha(\cdot)$ is increasing we have $s_*<s_0$.
On the other hand, in view of condition (i) in the statement of the Lemma and the finiteness of the second moment we have $s_0\alpha\to 0$. Thus,
\begin{equation}
s_*^+\alpha_{L,N}\to 0.
\label{s+}
\end{equation}
If $s_*<0$ we still have
\[
0\le \sum_k\big(k-\frac{N}{L}\big)\big(e^{s_*N/L}-e^{s_*k}\big)p_k=e^{s_*N/L}\big(\rho_c-\frac{N}{L}\big)\le \big(\rho_c-\frac{N}{L}\big)\longrightarrow 0,
\]
and since all the terms in the preceding sum are non negative, this implies
\begin{equation}
s_*^-\to 0.
\label{s-}
\end{equation}
The limits in (\ref{s+}), (\ref{s-}), together with the dominated convergence theorem and Fatou's lemma imply that
\begin{equation}
Z_{\alpha}(s_*)=\sum_{k\le\alpha}e^{s_*k}p_k\longrightarrow 1\qquad\text{and}\qquad Z_{\beta}(s_*)=\sum_{k\le\beta}e^{s_*k}p_k\longrightarrow 1\qquad\text{as } L\to\infty.
\label{partlimit}
\end{equation}
This in turn gives after another application of the dominated convergence theorem that
\begin{equation}
\sigma_\alpha^2(s_*)=\frac{1}{Z_{\alpha}(s_*)}\sum_{k\le\alpha}\big(k-\frac{N}{L}\big)^2e^{s_*k}p_k\longrightarrow \sigma^2\qquad \text{and}\qquad\sigma_\beta^2(s_*)\longrightarrow \sigma^2\quad\text{as }L\to\infty .
\label{sigmalimit}
\end{equation}
\noindent
We proceed now with the proof of the assertion of the lemma. Given $x\in \R$, write
\begin{eqnarray*}
\P\left[M_L\le \b_L(x),\,S_L=N \right]=Z^L_{\b}(s_*)\,e^{-s_*
N}\,\P_{\b}(s_*)[S_L=N]\,,
\end{eqnarray*}
and
\begin{eqnarray*}
\P\left[M_L\le \a_{L,N},\,S_L=N\right]=Z^L_{\a}(s_*)\,e^{-s_*
N}\,\P_{\a}(s_*)[S_L=N]\,.
\end{eqnarray*}
By condition {\em (i)} in the statement of the lemma we have
\begin{align}
\mu_{L,N}\left[ M_L\le \b_L(x) \right]&\simeq\frac{\P\left[M_L\le \b_L(x),\, S_L=N \right]}{\P\left[M_L\le\a_{L,N},\,S_L=N \right]}\notag\\
&=\left( \frac{Z_{\b}(s_*)}{Z_{\a}(s_*)}\right)^L \frac{\P_{\b}(s_*)[S_L=N]}{\P_{\a}(s_*)[S_L=N]}\,.
\label{limit}
\end{align}
By the local limit theorem for triangular arrays (Theorem 1.2 in \cite{DMcD}) and (\ref{sstar}) and (\ref{sigmalimit}), we have
\begin{equation}
\sqrt{2\pi L\sigma^2}\,\,\P_{\a}(s_*)\big[S_L=N\big]\longrightarrow 1\,.
\label{llt}
\end{equation}
In order to compute the asymptotics of $\P_{\b}(s_*)\big[S_L=N\big]$ in (\ref{limit}) we need to obtain estimates on $\rho_{\b}(s_*)-N/L$.
\begin{align*}
\rho_{\b}(s_*)=\frac{\sum_{k\le \b} ke^{s_*k}p_k }{\sum_{k\le \b} e^{s_* k}p_k} &=
\frac{\frac{N}{L}\sum_{k\le \a}e^{s_* k}p_k\,-\,\sum_{\b< k\le \a}ke^{s_* k}p_k }{\sum_{k\le \b}e^{s_*k}p_k}\\
&=\frac{N}{L}\,+\,\frac{\sum_{\b< k\le \a}\left(\frac{N}{L}-k\right)e^{s_* k}p_k }{\sum_{k\le \b}e^{s_*k}p_k}
\end{align*}
and
\begin{align}
L\left( \rho_{\b}(s_*)-\frac{N}{L}\right)^2 &=L\left(\frac{\sum_{\b< k\le \a}\left(\frac{N}{L}-k\right)e^{s_* k}p_k }{Z_{\b}(s_*)} \right)^2\notag\\
&\le \frac{L\,\,\sum_{\b< k\le \a} e^{s_* k} p_k}{Z^2_{\b}(s_*)}\, 
\sum_{\b<k\le \a}\left(\frac{N}{L}-k\right)^2e^{s_* k}p_k.
\label{density}
\end{align}
It now follows easily from (\ref{partlimit}), (\ref{sigmalimit}) and condition (ii) of the Lemma that
%
\[
\lim_{L\to \infty} L\left( \rho_{\b}(s_*)-\frac{N}{L}\right)^2\,=\,0\,.
\]
By another application of the local limit theorem for triangular arrays, we get that
\[
\sqrt{2\pi L\sigma^2}\,\,\P_{\b}(s_*)[S_L=N] \longrightarrow 1\qquad \text{as}\quad L\to \infty\,,
\]
and using condition (ii) (\ref{limit}) becomes
\[
\mu_{L,N}[M_L\le \beta_L(x)]\simeq \left(\frac{Z_{\b}(s_*)}{Z_{\a}(s_*)}\right)^L
=\left(1-\frac{\sum_{\b< k\le \a}e^{s_* k} p_k}{\sum_{k\le \a}e^{s_* k} p_k}\right)^L 
\longrightarrow e^{-\Phi(x)}\, .
\]
\CQFD
\subsection{The power law case}

\subsubsection{Proof of Theorem \ref{th1}}

Here $N>\rho_cL$ and $N/L\to \rho_c$ as $L\uparrow \infty$.\\

\noindent We will use that if $\lambda =1$ then $p_n\simeq A_1 n^{-b}$, and the decomposition (\ref{gammal1})
\begin{equation*}
N=\rho_c L+\sigma\sqrt{(b-3)L\log
L}\,\Big(1+\frac{b}{2(b-3)}\frac{\log\log L}{\log
L}+\frac{\gamma_L}{\log L}\Big)
\end{equation*}
that holds for $b>3$.\\

\noindent The proof of the theorem relies on (\ref{split}) and the following two lemmas. 
\begin{lemma}
Let $b>3$ and $N>\rho_c L$  be such that the sequence $\gamma_L$ in (\ref{gammal1}) has a limit
$\displaystyle{\lim_{L\to \infty} \gamma_L=\gamma \in [-\infty, \infty)\,}.$ If $\a_L=\frac{\sqrt{L}}{\log L}$ then
\[
\P\big[M_L\le \alpha_L; S_L=N \big]=\frac{1}{\sqrt{2\pi \sigma^2 L}}
\,\exp\left\{-\frac{(N-\rho_cL)^2}{2\sigma^2L} \right\}(1+o(1))\,.
\]
%
%
%
\label{lemma1}
\end{lemma}
\begin{lemma}
Suppose $b>3$ and $N-\rho_c L\gg\vartheta_L\sqrt{L}\,$, for a sequence $\vartheta_L\to\infty$. Then as $L\to\infty$
\[
\P\Big[ M_L\ge N-\rho_cL-\vartheta_L\sqrt{L};\ S_L=N\Big]=A_1 L(N-\rho_cL)^{-b} \big(1+o(1)\big).
\]
\label{lemma2}
\end{lemma}
\vspace{-.4cm}
{\em Proof of Lemma \ref{lemma1}:} The argument follows the standard approach used for moderate deviations of the sum of i.i.d. random variables. Consider $s_*\in \R$ such that $\rho_{\a_L}(s_*)=\frac{N}{L}$. Notice that $\rho_{\a_L}(0)=\rho_c-Z_{\a_L}^{-1}(0)\sum_{k>\a_L}(k-\rho_c)p_k<\frac{N}{L}$ for sufficiently large $L$, and by (\ref{s+}) we must have $s_*=o(\log L/\sqrt{L}).$ In particular, we have $Ls_*^{2+\epsilon}\to 0$ for all $\epsilon>0$.  Just as in the proof of Lemma \ref{lemma0}, we may write
\begin{align}
\P\left[ M_L\le \a_{L};\,S_L=N \right]&=Z^L_{\a_L}(s_*)\,e^{-s_* N}\,\P_{\a_L}(s_*)[S_L=N]\,\notag\\
&=Z_{\a_L}^L(0)\,\exp\left\{-L\int_0^{s_*} t\,\rho'_{\a_L}(t)\,dt \right\}\,\P_{\a_L}(s_*)[S_L=N]\,,
\label{ch.m}
\end{align}
\noindent
where we have used the identity
\[
\log\Big(\frac{Z_{\a_L}(s_*)}{Z_{\a_L}(0)}\Big)=\int_0^{s_*}\rho_{\a_L}(t)\,dt=s_*\rho_{\a_L}(s_*)-\int_0^{s_*}t\rho'_{\a_L}(t)\,dt\,.
\]
We now determine the asymptotic order of each term in (\ref{ch.m}). Observe that
\begin{equation*}
Z_{\a_L}^L(0)=\nu_{\phi_c}\left[ M_L \le \a_L \right]=\left(1-\sum_{k>\a_L}p_k\right)^L\longrightarrow 1\qquad\text{as }L\to \infty\,.
\end{equation*}
Furthermore, if we define $h_{\a_L}(t)=\rho_{\a_L}(t)-\rho_{\a_L}(0)-t\sigma_{\a_L}^2(0)$ we have
\begin{align}
\int_0^{s_*} t\rho_{\a_L}'(t)\,dt &=\frac{1}{\sigma_{\a_L}^2(0)}\int_0^{s_*}(\rho_{a_L}(t)-\rho_{\a_L}(0)-h_{\a_L}(t))\rho'_{\a_L}(t)\ dt\notag\\
&=\frac{\big(\rho_{\a_L}(s_*)-\rho_{\a_L}(0)\big)^2}{2\sigma_{\a_L}^2(0)}-\frac{1}{\sigma_{\a_L}^2(0)}\int_0^{s_*}h_{\a_L}(t)\rho_{\a_L}'(t)\ dt\,.
\label{h1}
\end{align}
Using elementary estimates one can show that
\begin{equation}
\frac{L\big(\rho_{\a_L}(s_*)-\rho_{\a_L}(0)\big)^2}{2\sigma^2_{\a_L}(0)}-\frac{(N-\rho_cL)^2}{2\sigma^2 L}\longrightarrow 0\,.
\label{h2}
\end{equation}
For the rightmost term in (\ref{h1}) notice that for all $s\in [0,s_*]$ we have
\[
0\le\sum_{k\le\a_L}k^2\big(e^{sk}-1\big)p_k\le Cs^\epsilon,
\]
for some $\epsilon>0$, which implies that $\big|\sigma_{a_L}^2(s)-\sigma_{\a_L}^2(0)\big|\le Cs^\epsilon$. Therefore,
\begin{align*}
&L\int_0^{s_*}h_{\a_L}(t)\rho_{\a_L}'(t)\ dt=L\int_0^{s_*}\left(\int_0^{t}\sigma_{\a_L}^2(s)-\sigma_{\a_L}^2(0)\ ds\right)\sigma_{\a_L}^2(t)dt=O\big(Ls_*^{2+\epsilon}\big)\to 0,
\end{align*}
where in the first equality we used that $\rho'_{\a_L}(t)=\sigma^2_{\a_L}(t)$. Together with (\ref{h1}), (\ref{h2}) this gives
\begin{equation*}
\exp\left\{-L\int_0^{s_*} t\rho'_{\a_L}(t)\ dt\right\}\simeq \exp\left\{-\frac{(N-\rho_cL)^2}{2\sigma^2 L}\right\}\qquad \text{as } L\to\infty\,.
\label{middle}
\end{equation*}
The assertion now follows recalling that $\sqrt{2\pi L\sigma^2_{\a_L}(s_*)}\P_{\a_L}(s_*)\big[S_L=N\big]\longrightarrow 1$
by (\ref{llt}).
\CQFD
\\

\noindent {\em Proof of Lemma \ref{lemma2}:} 
Consider a sequence $\vartheta_L$ as in the statement of the lemma. Then,
\begin{align*}
P_L& :=\P \big[ M_L \geq N-\rho_c L-\vartheta_L\sqrt{L};\, S_L =N\big] \\
& \simeq\sum_{k\ge(N-\rho_c L)-\vartheta_L\sqrt{L}} L\, p_k\,\P\big[S_{L-1} =N-k; \, M_{L-1}\leq k\big]\nonumber.
\end{align*}
Using the central limit theorem we can see that the contribution to the sum of the terms outside the set $ U_L=\{k\in\Z:\ |N-\rho_c L-k|\le \vartheta_L\sqrt{L}\}$ is negligible, that is
\[
P_L=\sum_{k\in U_L} L\, p_k\,\P\big[S_{L-1} =N-k; \, M_{L-1}\leq k\big] +o \big(L(N-\rho_cL)^{-b}\big)\,.
\]
We can now use the regular variation of $p_k$ to get
\[
P_L= A_1 L(N-\rho_cL)^{-b}\bigg(\sum_{k\in U_L} \P\big[S_{L-1} =N-k; \, M_{L-1}\leq k\big] +o(1)\bigg).
\]
The last sum converges to $1$ again by the central limit theorem and the fact that
$\P\Big[M_{L-1}\le k\Big]\to 1$, uniformly for $k\in U_L$, so
\[
P_L=A_1 L(N-\rho_cL)^{-b}\big(1+o(1)\big),
\]
as asserted.
\CQFD
\\

%
%
\noindent
We proceed now with the {\em proof of Theorem \ref{th1}}:\\

\noindent {\em a) The case} $\gamma_L\to-\infty$.\\

\noindent By Lemma \ref{lemma1} and (\ref{split}) if $N-\rho_cL\gg\sqrt{L}$ or the local limit theorem otherwise, condition {\em i}) of Lemma \ref{lemma0} is satisfied by
$\a_L=\sqrt{L}/\log L$. Consider $s_*>0$ such that $\rho_{\a_L}(s_*)=N/L$ and let $\b_L(x)=xL^{\frac{1}{L-1}},\,  x>0\,.$ Then
\begin{align*}
L\sum_{\b_L(x)<k\le \a_L} e^{s_*k}p_k&=L\sum_{\b_L(x)<k\le \a_L} p_k+L\sum_{\b_L(x)<k\le a_L}\big(e^{s_*k} -1 \big)\,p_k \\
&=L\bar{F}\Big(xL^{\frac{1}{b-1}}\Big)+L\sum_{\b_L(x)<k\le \a_L}\big(e^{s_*k} -1 \big)\,p_k - L\sum_{k>\a_L}p_k.\\
\end{align*}
It is easy to see that the first term above converges to $\frac{A_1}{b-1} x^{1-b}$ and that the last two terms vanish in the limit, since $s_*\a_L\to 0$ by (\ref{s+}). That is
\[
L\sum_{\b_L(x)<k\le a_L} e^{s_*k}p_k\longrightarrow \frac{A_1}{b-1} x^{1-b}\qquad \text{as }L\to \infty\, ,
\]
which is condition {\em ii}) in Lemma \ref{lemma0}.\\

\noindent
{\em b) The case} $\gamma_L\to+\infty$. \\

\noindent This case is essentially treated in \cite{armendarizetal08}. It is shown there (cf. Theorem 1b) that when the second term on the right hand side of (\ref{split}) dominates the probability of the event $\{S_L=N\}$, the variables $\{\eta_i\}$ aside from their maximum become asymptotically independent with distribution $\nu_{\phi_c}$. This entails that for all $y\in\R$
\[
\mu_{L,N}\Big[\frac{M_L -(N-\rho_c L)}{\sigma\sqrt{L}}\leq y\, \Big]\longrightarrow\frac{1}{\sqrt{2\pi }}\int_{-\infty}^y e^{-x^2 /2}\ dx,
\]
which is (\ref{fluct}). \\
\\
\noindent {\em c) The case} $\gamma_L\to\gamma\in\R$.\\

%
\noindent Here $N-\rho_c L\simeq\sigma\sqrt{(b-3)L\log L}\,$, and the two terms in the right hand side of (\ref{split}) are of the same order. Precisely,
\begin{equation}\label{lgdef}
\frac{\frac{1}{\sqrt{2\pi L\sigma^2}}\exp\big\{-\frac{(N-\rho_c L)^2}{2\sigma^2 L}\big\}}{LA_1 [N-\rho_cL]^{-b}}\longrightarrow\frac{\sigma^{b-1}(b-3)^{\frac{b}{2}}}{\sqrt{2\pi}A_1 }e^{-(b-3)\gamma}=:\ell_\gamma.
\end{equation}
It follows by (\ref{split}) and Lemma \ref{lemma1} that
\begin{equation*}
\liminf_{L\to \infty} \mu_{L,N}\Big[M_L\le \a_L\Big]\ge \frac{\ell_\gamma}{1+\ell_\gamma}.
\end{equation*}
On the other hand, applying Lemma \ref{lemma2} with $\vartheta_L \ll\sqrt{\log L}$ we have that
\begin{equation*}
\liminf_{L\to \infty}\mu_{L,N}\Big[M_L\ge N-\rho_cL-\vartheta_L\sqrt{L}\Big]= \frac{1}{1+\ell_\gamma}\,,
\end{equation*}
and by the central limit theorem
\[
\limsup_{L\to \infty}\mu_{L,N}\Big[M_L\ge N-\rho_cL+\vartheta_L\sqrt{L}\Big]= 0\,.
\]
The last three relations together imply that
\begin{equation}
\frac{M_L}{N-\rho_c L} \stackrel{\mu_{L,N}}{\Longrightarrow} Be(p_\gamma)\quad \text{with }\ p_\gamma=\frac{1}{1+\ell_\gamma}.
\label{bernoulimit}
\end{equation}

\noindent
Note that $p_\gamma$ as given above satisfies
\[
p_\gamma \to\left\{\bacl 1\ &\mbox{ if }\gamma\to\infty\\ 0\ &\mbox{ if }\gamma\to -\infty\ea\right.\ .
\]
This finishes the proof of Theorem \ref{th1}.
\CQFD

\subsubsection{Proof of Theorem \ref{th2}}

\noindent Here $N<\rho_c L$ and $N/L \to \rho_c$ as $L\uparrow \infty$.\\

\noindent We take $\a_L=\infty$, so that ({\em i}) in Lemma \ref{lemma0} is automatically satisfied. 
It remains to identify
the sequence $\b_L\to \infty$ and the limit $\Phi$ in condition ({\em ii}) for
each case
({\em a}), ({\em b}) or ({\em c}) in the theorem. Note that since
\[
\frac{N}{L}=\rho(s_*)=\rho_c+\int_0^{s_*}\sigma^2(s)\ ds, \qquad\text{with}\qquad \sigma^2(s)\stackrel{s\to 0}{\longrightarrow}\sigma^2,
\]
we have
\begin{equation}
s_*=\frac{N-\rho_cL}{\sigma^2L}\big(1+o(1)\big)\qquad \text{as } L\to\infty.
\label{sas}
\end{equation}
\noindent{\em a)} The case $\omega_L\to 0$.\\

\noindent Let $\displaystyle{\b_L(x)=xL^{\frac{1}{b-1}}\,.}$ Then
\[
L\sum_{k>\b_L(x)}e^{s_*k} p_k=L\bar{F}\big[\b_L(x)\big]+L\sum_{k>\b_L(x)}\big(e^{s_*k}-1 \big)p_k\,,
\]
where
\[
0\ge L\sum_{k>\b_L(x)}\big(e^{s_*k}-1 \big)p_k\ge Ls_*\sum_{k>\b_L(x)}kp_k =O(\omega_L)\longrightarrow 0\,.
\]
That is,
\[
L\sum_{k>\b_L(x)}e^{s_*k} p_k \longrightarrow
\frac{A_1}{b-1}\, x^{1-b}\,.
\]

\noindent{\em b)} The case $\omega_L\to \omega >0$.\\

\noindent As in the previous case, let
$\displaystyle{\b_L(x)=xL^{\frac{1}{b-1}}}\,$. By the regular
variation of the probabilities $p_k$,
\begin{equation*}
L\sum_{k>\b_L(x)} e^{s_*k}p_k\,\simeq\,
A_1 L\sum_{k=\b_L(x)+1}^{ML^{\frac{1}{b-1}}} e^{s_*k}\frac{1}{k^b}\,+\,
A_1 L\sum_{k>ML^{\frac{1}{b-1}}}e^{s_*k}\frac{1}{k^b}\,.
\end{equation*}
We compute the limits of both terms on the right hand side above:
\[
\lim_{M\to \infty}\lim_{L\to
\infty}L\sum_{k>ML^{\frac{1}{b-1}}}e^{s_*k}\frac{1}{k^b}\,=\,0
\]
and by (\ref{sas})
\begin{align*}
\lim_{M\to \infty}\lim_{L\to \infty}L\sum_{k=\b_L(x)+1}^{ML^{\frac{1}{b-1}}}
e^{s_*k}\frac{1}{k^b}\,&=\,\lim_{M\to \infty}\lim_{L\to \infty}
L\sum_{k=\b_L(x)+1}^{ML^{\frac{1}{b-1}}} \exp\left\{-\omega\frac{k}{L^{\frac{1}{b-1}}}\right\} \frac{1}{k^{b}}\\
&=\lim_{M\to \infty} \int_{x}^{M}e^{-\omega t} \frac{1}{t^b}\,dt\,.
\end{align*}

\noindent{\em c)} The case $\omega_L\to \infty$.\\

\noindent Define now a sequence $B_L$ by the equation
\begin{equation}
\left(|s_*|B_L\right)^b e^{|s_*|B_L}=A_1 L|s_*|^{b-1}\,,
\label{BsubL}\end{equation} and note that $\omega_L\to \infty$
implies that $|s_*|B_L\to \infty$ as well. Let
$\displaystyle{\b_L(x)=B_L+\frac{x}{|s_*|}}\,.$ Then
\begin{align*}
L\sum_{k>\b_L(x)} e^{s_* k}p_k&=Le^{s_*
B_L}\sum_{k>\b_L(x)}e^{s_*(k-B_L)} p_k \\
&\simeq A_1 Le^{-|s_*| B_L}\sum_{k>\b_L(x)}e^{s_*(k-B_L)}\frac{1}{k^b}\\
&=|s_*|\sum_{k>\b_L(x)}e^{s_*(k-B_L)}\left(\frac{B_L}{k}\right)^b\\
&=|s_*|\sum_{k'>x/|s_*|}e^{s_*(k')}\frac{|s_*|B_L}{|s_*|k' +|s_*|B_L}\longrightarrow \int_x^{\infty}
e^{-t}\,dt=e^{-x}
\end{align*}
as $L\to\infty$, using dominated convergence with $|s_*|B_L\to\infty$. The third line above follows from the second one by (\ref{BsubL}).
This concludes the proof of the theorem.
 \CQFD

\subsection{The stretched exponential case} Here we have $p_n\simeq A_\lambda e^{-\frac{b}{1-\lambda}n^{1-\lambda}}$.
The proofs in this case use results from  \cite{nagaev68}, which are summarized in the Appendix. \\

\subsubsection{Proof of Theorem \ref{th3}}

We recall the the notation from equation (\ref{gammal2})
\[
N=\rho_cL+t_L(\sigma^2L)^{\frac{1}{1+\lambda}}\,.
\]
{\em a) The case} $\displaystyle{\limsup_{L\to \infty} t_L<c_\lambda}$. \\

\noindent The second equation in (\ref{mdal}) that gives the limit theorem for $M_L$ without conditioning
is a standard computation in extreme value theory. The appropriate scales
\[
y_L\simeq \left( \frac{1-\lambda}{b} \log L\right)^\frac{1}{1-\lambda},\qquad b_L\simeq \frac{y_L^{\lambda}}{b}\qquad\text{as }\quad L\to \infty
\]
are chosen so that
\begin{equation}\label{ybscale}
 L \sum_{k > y_{L}+x b_L}p_k \to e^{-x}, \quad x\in\R.
\end{equation}
Let $c>1$. According to (\ref{rem1})
\[
\P\big[M_L\le c y_L\ ; \ S_L=N \big]\simeq\P\big[ S_L=N\big]\,.
\]
We thus set the sequence $\a_L=c y_L$, and item {\em i})  in Lemma \ref{lemma0} is satisfied. Recall that $s_*>0$ since $N>\rho_cL$.
If we define $\b_L(x)=y_L+x b_L$, we obtain
\begin{align*}
\lim_{L\to \infty} L\sum_{\b_L(x)<k\le \a_L}e^{s_* k}p_k&=\lim_{L\to \infty} L\sum_{\b_L(x)<k\le \a_L} p_k \\
&=\lim_{L\to \infty} L\sum_{k>y_L+xb_L}p_k\,-\,\lim_{L\to \infty}L\sum_{k>c y_L}p_k \\
&=e^{-x}\,,
\end{align*}
where the first identity follows from (\ref{s+}) and the third one is (\ref{ybscale}). This provides condition {\em ii}) in Lemma
\ref{lemma0}.\\


\noindent {\em b) The case} $t_L \to t>c_\lambda$.\\

\noindent When $N-\rho_cL\gg L^{\frac{1}{2\lambda}}$ this can be deduced by Theorem 1a in \cite{armendarizetal08}, since in that case the $L-1$ smallest variables become asymptotically independent. In fact, we can then take $a_L=1$. For smaller values of $N-\rho_c L$, even though it is not stated explicitly, this is essentially proved in \cite{nagaev68}. Note that by Remarks 2,3 and 5 in the Appendix for any $\theta_L\to\infty$ we have
\[
\mu_{L,N}\Big[ |M_L-(1-\alpha)(N-\rho_cL)|<\theta_L\sqrt{L}\Big]\longrightarrow 1.
\]
where $\alpha$ is implicitly defined by (\ref{alef}). This implies that the conditional distribution of $\displaystyle\frac{M_L-(1-\alpha)(N-\rho_cL)}{\sqrt{L}}$
is tight. A careful reading of the proof of Lemma 7, part 2 in \cite{nagaev68} reveals that in fact this distribution has the asserted limit. Note that the sequence $a_L$ in the statement of Theorem \ref{th2} is given by $1-\alpha$. The properties of its limit $a(t)$ can easily be deduced from (\ref{alef}).\\
\\
{\em (c) The case} $t_L \to c_\lambda$. \\

\noindent This case can be treated analogously to the third part of Theorem \ref{th1}. By (\ref{Pc4}), the leading order of $\P\big[S_L=N \big]$ is the sum of two explicit terms, and one has to find the precise subscale around $N-\rho_cL-c_\lambda(\sigma^2L)^{\frac{1}{1+\lambda}}$ where these two terms are of the same order. Using (\ref{simple}) for $\lambda >1/2$ and (\ref{alef}) we find that on the scale (\ref{subl})
\bea\label{pstretched}
\lefteqn{\frac{1}{\sigma\sqrt{2\pi L}}\, e^{-\frac{k^2}{2L\sigma^2}}\bigg/\frac{A_\lambda L}{\sqrt{1-\frac{\sigma^2\gamma\l(1-\l)L}{k^{1+\l}(1-\alpha)^{1+\l}}}}\, e^{-\frac{\alpha^2k^2}{2\sigma^2L}-\gamma(1-\alpha)^{1-\l}k^{1-\l}}}\nonumber\\
& &\quad\longrightarrow \frac{\sqrt{1+\lambda}}{2A_\lambda \sqrt{\pi\sigma^2}}\, e^{\gamma }=:\ell_\gamma\ .
\eea
From this, (\ref{bern2}) can be deduced analogous to the proof of Theorem \ref{th1} with $p_\gamma =(1+\ell_\gamma )^{-1}$.
\CQFD

\subsubsection{Proof of Theorem \ref{th4} }

As for the downside deviations in the power law case, we here set $\a_L=\infty$. Then {\em i}) in Lemma \ref{lemma0} is satisfied, $s_*<0$,
and satisfies (\ref{sas}).
%
\noindent Now the limit  in {\em ii}), Lemma \ref{lemma0} is  given by $\Phi(x)=e^{-x}$ in all cases, so we just need to prove that the proposed values of the sequence $\b_L(x)$ do the job. The proofs are all based on the following computation, with simple adjustments to match each situation.\\

\noindent Let $\gamma_L$ and $\zeta_L$ be sequences such that there exist $\ell_1$ and $\ell_2\,\in \R$ with 
\bea
\lim_{L\to \infty} \frac{\zeta_L^2}{\gamma_L^{1+\lambda}}=0,\quad \lim_{L\to \infty} \zeta_L |s_*|=\ell_1\,\quad \text{and } \quad
\lim_{L\to \infty} \frac{\zeta_L}{\gamma_L^{\lambda}}=\ell_2,\quad \ell_1+b\ell_2=1.
\label{gamma-zeta}
\eea
Then
\begin{eqnarray*}
\lefteqn{A_\lambda L\sum_{k= \gamma_L+x\zeta_L}^{\gamma_L+y\zeta_L} e^{s_* k}e^{-\frac{b}{1-\lambda}k^{1-\lambda}}=
A_\lambda L\,e^{s_*\gamma_L-\frac{b}{1-\lambda}\gamma_L^{1-\lambda} }
\sum_{k= \gamma_L+x\zeta_L}^{ \gamma_L+y\zeta_L} e^{s_* (k-\gamma_L)}e^{-\frac{b}{1-\lambda}\big(k^{1-\lambda}-\gamma_L^{1-\lambda}\big)} }\notag\\
& &\qquad =A_\lambda L\,e^{s_*\gamma_L-\frac{b}{1-\lambda}\gamma_L^{1-\lambda} }
\sum_{k= \gamma_L+x\zeta_L}^{\gamma_L+y\zeta_L} e^{s_* (k-\gamma_L)}e^{-\frac{b}{\gamma_L^{\lambda}}(k-\gamma_L)}\,\big(1+o(1)\big)\notag\\
& &\qquad =A_\lambda L\,e^{s_*\gamma_L-\frac{b}{1-\lambda}\gamma_L^{1-\lambda} }
\sum_{k-\gamma_L=x\zeta_L}^{k-\gamma_L=y\zeta_L}
e^{-(|s_*|+\frac{b}{\gamma_L}) (k-\gamma_L)}\,\big(1+o(1)\big)\notag\\
& &\qquad\simeq A_\lambda L\,\frac{e^{s_*\gamma_L-\frac{b}{1-\lambda}\gamma_L^{1-\lambda} }}{|s_*|+\frac{b}{\gamma_L^{\lambda}}}
\int_{x(\ell_1+b\ell_2)}^{y(\ell_1+b\ell_2)} e^{-t}\,dt \,.
\end{eqnarray*}
We now let $y\to \infty$ and apply (\ref{gamma-zeta}) to get
\begin{align}
A_\lambda L\sum_{k> \gamma_L+x\zeta_L} e^{s_* k}e^{-\frac{b}{1-\lambda}k^{1-\lambda}}
&\simeq A_\lambda L\,\,\gamma_L^{\lambda}\,\,\frac{e^{s_*\gamma_L-\frac{b}{1-\lambda}\gamma_L^{1-\lambda} }}{|s_*|\gamma_L^{\lambda}+ b}\, \,e^{-x}\,. 
 \label{exponential-phi}
\end{align}
We will work on this limit case by case.\\

\noindent Recall from the proof of Theorem \ref{th3} the existence of sequences
\[
y_L\simeq \left( \frac{1-\lambda}{b} \log L\right)^\frac{1}{1-\lambda},\qquad b_L\simeq \frac{y_L^{\lambda}}{b}\qquad\text{as }\quad L\to \infty
\]
 satisfying 
 \bea
 L \sum_{k > y_{L}+x b_L}p_k \to e^{-x}, \,\,\, x\in\R.
\label{ybscale2}
 \eea

\noindent{\em a) The case $s_*y_L\to 0\Leftrightarrow \omega_L\to 0$.} \\

\noindent Here $\gamma_L:=y_L$ and $\zeta_L:=b_L$ satisfy (\ref{gamma-zeta}) with $\ell_1=0$
and $\ell_2=1/b$. In fact, in this case it is straightforward from (\ref{ybscale2}) and dominated convergence that
\[
\lim_{L\to \infty} L\sum_{k>y_L+xb_L}e^{s_* k}p_k=e^{-x}.
\]
\noindent{\em b) The case $s_* y_L^{\lambda}\to 0$.}\\

\noindent  Here we need to be slightly more careful with the choice of $\gamma_L$. Set $\zeta_L=b_L$ and let $\gamma_L$ be the solution to
\bea
A_\lambda L\,\,b_L\,\,e^{s_*\gamma_L-\frac{b}{1-\lambda}\gamma_L^{1-\lambda} }=1
\label{case2}
\eea
The sequences $y_L$ and $b_L$ can be chosen so that
\[
A_\lambda L\,b_Le^{-\frac{b}{1-\lambda}y_L^{1-\lambda}}=1,
\]
and from (\ref{case2})
\[
s_*\gamma_L=\frac{b}{1-\lambda}\big(\gamma_L^{1-\lambda}-y_L^{1-\lambda}\big)\qquad\text{or}\qquad s_*y_L^{\lambda}=
\frac{b}{1-\lambda}\left(\frac{y_L^{\lambda}}{\gamma_L^{\lambda}}-\frac{y_L}{\gamma_L} \right)\,.
\]
By the condition $s_*y_L^{\lambda}\to 0$ this implies that $\displaystyle{\frac{y_L}{\gamma_L}\simeq 1}$, and (\ref{gamma-zeta}) holds with
$\ell_1=0$ and $\ell_2=1/b$. Also (\ref{case2}) and (\ref{exponential-phi}) imply 
\[
\lim_{L\to \infty} L \sum_{k>\gamma_L+xb_L} e^{s_* k} p_k=e^{-x}\,.
\]

\noindent {\em c) The case $s_* y_L^{\lambda} \to c<0,\, c\in \R$. }\\

\noindent The scaling in the sequences $\gamma_L$ and $\zeta_L$ is preserved, but the limits $\displaystyle{\lim_{L\to \infty} \frac{\gamma_L}{y_L}}$ and $\displaystyle{\lim_{L\to \infty}\frac{\zeta_L}{b_L}}$ need to be chosen so that $\ell_1$ and $\ell_2$ in (\ref{gamma-zeta})
satisfy $\ell_1+b\ell_2=1$ and the right hand side of (\ref{exponential-phi}) equals $e^{-x}\,$,
\[
A_\lambda L\,\,\gamma_L^{\lambda}\,\,\frac{e^{s_*\gamma_L-\frac{b}{1-\lambda}\gamma_L^{1-\lambda} }}{|s_*|\gamma_L^{\lambda}+ b}=1\,
\]
{\em d) The case $|s_*| y_L^{\lambda}\to \infty$.}\\
\noindent Let now $\displaystyle{\zeta_L=\frac{1}{|s_*|}}$ and set $\gamma_L$ as the solution to 
\bea
\frac{A_\lambda L}{|s_*|}\,\,e^{s_*\gamma_L-\frac{b}{1-\lambda}\gamma_L^{1-\lambda} }=1\,.
\label{last}
\eea
(It is easy to see that such a solution exists). Now taking logarithms in (\ref{last}) we obtain that, to leading order,
\[
\left(\frac{1-\lambda }{b} \log L\right)^{\frac{1}{1-\lambda}} \simeq
\gamma_L\left(1+\frac{1-\lambda}{b}|s_*|\gamma_L^{\lambda}\right)^{\frac{1}{1-\lambda}},
\]
from where we conclude that necessarily $|s_*|\gamma_L^{\lambda}\to \infty$. Then (\ref{gamma-zeta}) holds with $\ell_1=1$ and $\ell_2=0$, and (\ref{exponential-phi}) follows from (\ref{last}).
\CQFD
\subsection{Proof of Theorem \ref{bulk}}
a) The assertion will follow from Theorem 24.2 in \cite{bill} provided we check the validity of the following three conditions for the exchangeable random variables $\xi_{x}^{L,N}= \frac{\eta_x-N/L}{\sigma\sqrt{L}}$.\\
1. $\displaystyle\sum_{x=1}^L \xi_{x}^{L,N}\stackrel{\mu_{L,N}}{\longrightarrow}0.$ This is trivial since the sum of all $\xi_x^{L,N}$ is equal to zero $\mu_{L,N}$-a.s.\\
2. $\displaystyle \big|\max_{1\le x\le L}\xi_{x}^{L,N}\big|\stackrel{\mu_{L,N}}{\longrightarrow}0.$ This follows from part (a) of Theorems \ref{th1}, \ref{th3} when $N\ge\rho_c L$, or from Theorems \ref{th2}, \ref{th4} when $N<\rho_c L$.\\
3. $\displaystyle \sum_{x=1}^L \big(\xi_{x}^{L,N}\big)^2\stackrel{\mu_{L,N}}{\longrightarrow} 1$: We prove it in detail for the case when the occupation variables $\eta_x$ follow a power law,  the stretched exponential case being completely similar.

Let $\epsilon>0$ and set $\alpha_L=\frac{\sqrt{L}}{\log L}$ as in Lemma \ref{lemma1} . Then
\begin{align}
R_{L,N}&:= \mu_{L,N}\left[\Big|\frac{1}{L}\sum_{x=1}^L\big(\eta_x-\frac{N}{L}\big)^2-\sigma^2 \Big|>\epsilon \right] \notag\\
& \le \mu_{L,N} \left[\Big| \frac{1}{L}\sum_{x=1}^L\big(\eta_x-\frac{N}{L}\big)^2-\sigma^2 \Big|>\epsilon,\, M_L \le \alpha_L \right] 
+\mu_{L,N}\big[ M_L>\alpha_L\big]
\label{a1}
\end{align}
By Theorem \ref{th1} a) and Theorem \ref{th2} the second term on the right side above tends to $0$ as $L\to \infty$. Let us now write
\begin{align}
&\mu_{L,N}\left[\Big|\frac{1}{L}\sum_{x=1}^L\big(\eta_x-\frac{N}{L}\big)^2-\sigma^2 \Big|>\epsilon,\, M_L\le \alpha_L \right] \notag\\
&\qquad=\frac{\mu^L\left[\big|\frac{1}{L}\sum_{x=1}^L \big(\eta_x- \frac{N}{L}\big)^2-\sigma^2\big|> \epsilon,\,M_L\le \alpha_L,\,
\sum_{x=1}^L \eta_x=N \right]}{\mu^L\left[ \sum_{x=1}^L \eta_x= N\right]}\notag\\
&\qquad\qquad=\frac{\P_{\alpha_L}(s_*)\left[\big|\frac{1}{L}\sum_{x=1}^L \big(\eta_x- \frac{N}{L}\big)^2-\sigma^2\big|> \epsilon,\,
\sum_{x=1}^L \eta_x=N \right]}{ Z_{\alpha_L}^L(s_*)\,\E\left[e^{s_*S_L}\,\mathbbm{1}_{\{\sum \eta_x =N\}}\right]} \notag\\
&\qquad\qquad\qquad\le \frac{\P_{\alpha_L}(s_*)\left[\big|\frac{1}{L}\sum_{x=1}^L \big(\eta_x- \frac{N}{L}\big)^2-\sigma^2\big|> \epsilon 
\right]}{\P_{\alpha_L}(s_*)\left[\sum_{x=1}^L \eta_x=N\right]}\,,
\label{a2}
\end{align}
where we recall that given parameters $\alpha>0$ and $s\in \R$, the measure $\P_{\alpha}(s)$ is defined by (\ref{tilted}).
From (\ref{a1}), (\ref{a2}) and (\ref{llt}) we conclude that
\[
R_{L,N}
\le \sqrt{2\pi \sigma^2 L} \,\P_{\alpha_L}(s_*)\left[   \Big|\frac{1}{L}\sum_{x=1}^L \big(\eta_x- \frac{N}{L}\big)^2-\sigma^2\Big|> \epsilon \right]\,+o(1) \,.
\]
The result will thus follow if we show that
\bea
\sqrt{L}\,\,\P_{\alpha_L}(s_*)\left[ \frac{1}{L}\sum_{x=1}^L \big(\eta_x- \frac{N}{L}\big)^2 \ge \sigma^2 + \epsilon \right]\longrightarrow 0\quad
\mbox{ as }L\to \infty
\label{a3}
\eea
and
\bea
\sqrt{L}\,\,\P_{\alpha_L}(s_*)\left[ \frac{1}{L}\sum_{x=1}^L \big(\eta_x- \frac{N}{L}\big)^2\le \sigma^2 - \epsilon \right]\longrightarrow 0\quad
\mbox{ as }L\to \infty\,.
\label{a4}
\eea
Let us start with the former. If $\zeta>0$ then
\begin{align}
&\P_{\alpha_L}(s_*)\left[\sum_{x=1}^L \big(\eta_x- \frac{N}{L}\big)^2 \ge (\sigma^2 + \epsilon)L \right] \le e^{-\zeta (\sigma^2+\epsilon)L}\,
\E^{s_*}\left[ \exp\left\{ \zeta \sum_{x=1}^L \big(\eta_x-\frac{N}{L}\big)^2 \right\} \right]\notag\\
&\hspace{5cm}=e^{-\zeta (\sigma^2+\epsilon)L}\left(\frac{1}{Z_{\alpha_L}(s_*)}\sum_{k=1}^{\alpha_L} e^{\zeta(k-\frac{N}{L})^2} e^{s_*k}p_k\right)^L\,,
\label{a5}
\end{align}
where $\E^{s_*}$ denotes expectation with respect to the measure $\P_{\alpha_L}(s_*)$. \\
\\
Now set $\zeta=\frac{\log^{3/2} L}{L}$, so that
$\zeta \a_L^2 \to 0$, and apply the elementary inequality $e^x\le 1+x\psi(h)$  for $x\in[0,h]$, where $\psi(h)=\frac{e^h-1}{h}$. We get
\begin{align*}
\frac{1}{Z_{\alpha_L}(s_*)}\sum_{k\le \alpha_L} e^{\zeta (k-\frac{N}{L})^2}e^{s_*k}p_k
&\le \frac{1}{Z_{\alpha_L}(s_*)}\sum_{k\le \alpha_L} \left[1+\zeta \big(k-\frac{N}{L}\big)^2\psi(\zeta\a_L^2)\right] e^{s_*k}p_k\\
&=1+\zeta\, \sigma^2_{\alpha_L}(s_*)\psi(\zeta\a_L^2).
\end{align*}
From (\ref{sigmalimit}) we can make $\sigma^2_{\alpha_L}(s_*)\psi(\zeta\a_L^2)<\sigma^2+\epsilon/2$ for large enough $L$, so (\ref{a5}) becomes
\begin{align}
&\P_{\alpha_L}(s_*)\left[\sum_{x=1}^L \big(\eta_x- \frac{N}{L}\big)^2 \ge (\sigma^2 + \epsilon)L \right] \le e^{-\zeta (\sigma^2+\epsilon)L}
 \left(1+ \zeta\, \sigma^2_{\alpha_L}(s_*)\psi(\zeta\a_L^2)\big)  \right)^L \notag\\
 &\qquad\qquad\le e^{-\zeta (\sigma^2+\epsilon)L}\,\, e^{\zeta L\sigma^2_{\alpha_L}(s_*)\psi(\zeta\a_L^2)}\le e^{-\frac{\epsilon}{2}\log^{3/2}L},
\label{final estimate}
\end{align}
from where (\ref{a3}) is easily obtained. The limit (\ref{a4}) can be derived by similar estimates.\\

\noindent b) In the power law case ($\lambda=1$), or when $\lambda<1$ and $N-\rho_c L\gg L^{\frac{1}{2\lambda}}$, the assertion follows immediately (with $a_L=1$) from the asymptotic independence of the bulk variables proved in Theorems 1b, 1a in \cite{armendarizetal08}, respectively. \\
\\
Let us then consider the stretched exponential case $\lambda <1$ when $N-\rho_cL=t_L (\sigma^2 L)^{\frac{1}{1+\lambda}}$, and 
$t_L\to t \in (c_{\lambda},+\infty]$ (we refer to the statement of Theorem \ref{th3} for notation). The case 
$N-\rho_c L\gg L^{\frac{1}{2\lambda}}$ discussed in the previous paragraph clearly belongs to this family as well. \\
\\
We first observe that in this situation
\bea
N-a(t)(N-\rho_cL)<\rho_cL+c_{\lambda} (\sigma^2 L)^{\frac{1}{1+\lambda}}.
\label{reduction}
\eea
Indeed, according to (\ref{alef}), $a(t)$ satisfies
\[
\frac{1}{t^{1+\lambda}}=\frac{(1-a(t))\,a(t)^{\lambda}}{\gamma(1-\lambda)}\quad\mbox{with } \quad \gamma=\frac{b}{1-\lambda}\,.
\] 
Let $x_t=t(1-a(t))$. Then it follows from Theorem 2 in the Appendix that $x_t$ is the smallest positive root of the equation
\bea
b=x_t(t-x_t)^{\lambda},
\label{x_t}
\eea
and it is easily checked that 
\[
\lim_{t\uparrow \infty} x_t=0 \qquad \mbox{and} \qquad \lim_{t\downarrow c_{\lambda}} x_t=c_{\lambda} \,\frac{1-\lambda}{1+\lambda}<c_\lambda\,.
\]
In order to conclude (\ref{reduction}) it will therefore be enough to show that $x_t$ is decreasing. Differentiating in (\ref{x_t}) we get, after a couple of operations, that the derivative $x'_t$ satisfies
\[
x'_t \left(\frac{\lambda}{t-x_t}-\frac{1}{x_t} \right)=\frac{\lambda}{t-x_t}.
\]
Now 
\[
\frac{\lambda}{t-x_t}<\frac{1}{x_t} \iff \lambda x_t < t-x_t \iff \lambda t(1-a(t))<ta(t) \iff \frac{\lambda}{1+\lambda}<a(t)\,.
\]
But $a(t)$ is increasing on the half-line $(c_{\lambda}, \infty)$ with  $\lim_{t\downarrow c_{\lambda}} a(t)=\frac{2\lambda}{1+\lambda}$, $\lim_{t\uparrow \infty}a(t)=1$, and hence $a(t)>\frac{\lambda}{1+\lambda}$. \\
\\
Inequality (\ref{reduction}) allows us to decompose the random walk $Y_L$ into two components: the first term will be easily shown to converge to a Brownian bridge via the same arguments applied to prove the first statement in the theorem, while the second one is a drift term determined by the Gaussian limit specified in Theorem \ref{th3} b). Precisely, write
\bea
Y^L_s = W^L_s+\frac{[sL]}{L}\frac{M_L-a_L(N-\rho_cL)}{\sqrt{\sigma^2 L}}\,,\qquad W^L_s=\frac{1}{\sigma \sqrt{L}}\sum_{x=1}^{[sL]} \left( \tilde{\eta}_x-\frac{N-M_L}{L}\right).
\label{decompos}
\eea
Next, consider the interval  ${\cal A}_L=\left\{X\in \R,\,\Big|\frac{X}{a_L(N-\rho_cL)}-1\Big| \le \delta_L \right\}$ associated to $\delta_L=L^{-\frac{1}{4}\frac{1-\lambda}{1+\lambda}}$,  chosen 
 so that Theorem \ref{th3} b) implies $\lim_{L\to \infty} \mu_{L,N}[\,M_L\in {\cal A}_L]=1$. 
 Notice that  by (\ref{reduction}) the occupation variables $\{\tilde{\eta}_x\}_{x=1,\cdots,L}$ are in the subcritical regime when $M_L\in {\cal A}_L$, they are clearly exchangeable, and moreover when properly centered they satisfy conditions 1), 2) and 3) of Theorem 24.2 in \cite{bill}.\\
 \\
 Namely,
let $\tilde{\xi}_x^{L,N}=\frac{\tilde{\eta}_x-(N-M_L)/L}{\sigma\sqrt{L}}$. Then, provided $M_L \in {\cal A}_L$, we can easily show that\\
\\
1'. $\displaystyle{\sum_{x=1}^L \tilde{\xi}^{L,N}\stackrel{\mu_{L,N}}{\longrightarrow}0}$. In fact, the sum equals $0$ except in the rare event that the second order statistic of the sample $\{\eta_x\}_{1\le x\le L}$is greater than $L^{1/4}$.\\
\\
2'. $\displaystyle{|\max_{1\le x \le L}\tilde{\xi}_x^{L,N}|\stackrel{\mu_{L,N}}{\longrightarrow}0}$. This is trivial, as $\tilde{\eta}_x\le L^{1/4}$ and 
$N-M_L\le L^{\frac{1}{1+\lambda}}$ when $M_L\in {\cal A}_L$.\\
\\
3'. $\displaystyle{\sum_{x=1}^L (\tilde{\xi}^{L,N}_x)^2\stackrel{\mu_{L,N}}{\longrightarrow}1}$. This can be shown following the arguments applied to prove condition 3) in the first statement of the theorem, the difference being that it is now necessary to condition on $M_L$
before applying Chebyshev's inequality: for $\epsilon>0$,
\begin{align}
& \mu_{L,N}\left[\Big|\frac{1}{L}\sum_{x=1}^L\big(\tilde{\eta}_x-\frac{N-M_L}{L}\big)^2-\sigma^2 \Big|>\epsilon,\, M_L\in {\cal A}_L \right] \notag\\
& \hspace{1cm} = \E^{\mu_{L,N}} \left[\mu_{L,N}\left[\Big| \frac{1}{L}\sum_{x=1}^L\big(\tilde{\eta}_x-\frac{N-M_L}{L}\big)^2-\sigma^2 \Big|>\epsilon\Big|\,M_L\right]\, \mathbbm{1}_{\{M_L \in {\cal A}_L\}}\right]\,.
\label{condition-first}
\end{align}
The estimates leading to the bound (\ref{final estimate}) hold uniformly for $\tilde{N}=N-M_L$ when $M_L\in {\cal A}_L$, and hence can be applied to the conditioned expectation in the right side of (\ref{condition-first}) to conclude that 
\[
 \mu_{L,N}\left[\Big|\frac{1}{L}\sum_{x=1}^L\big(\tilde{\eta}_x-\frac{N-M_L}{L}\big)^2-\sigma^2 \Big|>\epsilon,\, M_L\in {\cal A}_L \right] \longrightarrow  0\,,
\]
as required. Conditions 1'), 2') and 3') imply that $W_s^L\stackrel{\mu_{L,N}}{\Longrightarrow} BB_s$, the standard Brownian bridge on $[0,1]$ conditioned to return to the origin at time $1$. On the other hand, we know form Theorem \ref{th3} b) that $\frac{M_L-a_L(N-\rho_cL)}{\sqrt{\sigma^2 L}}\stackrel{\mu_{L,N}}{\Longrightarrow} \Phi$, a zero mean Gaussian variable with variance $1/(1-\frac{\lambda(1-a(t))}{a(t)})$. \\
\\
Our assertion will follow once we prove the convergence of the finite dimensional marginals plus tightness for the laws of the sequence  $Y_s^L$.  The former is easily derived by first conditioning on $M_L$; the fact that the limit of the first term in (\ref{decompos}) is independent of the value of the second implies that the finite dimensional marginal distributions converge to those of $BB_s+s\,\Phi$.
Tightness is also straightforward: by the linearity of the second term in (\ref{decompos}), it suffices to show that the modulus of continuity of the first term tends to $0$,
\[
\omega(\delta)=\sup_{|s-r|\le \delta} \big| W^L_s-W^L_r\big|\,\stackrel{\mu_{L,N}}{\longrightarrow} 0 \qquad \mbox{as}\qquad\delta\to 0\,,
\]
which is a direct consequence of the Arzel\`a-Ascoli Theorem. \CQFD

\section*{Acknowledgments}
The results in this article answer a question raised by Pablo Ferrari after a seminar in the XI Brazilian School of Probability, and we would like to thank him for pointing out the problem to us. We are also grateful to Paul Chleboun who provided us with the simulation data for Figure~\ref{fig:profile}. This research has been supported by the University of Warwick Research Development Grant RD08138, the FP7-REGPOT-2009-1 project Archimedes Center for Modeling, Analysis and Computation, PICT-2008-0315 project "Probability and Stochastic
Processes" and the University of Crete Basic Research grant KA-2865. I.A.\ and S.G.\ are also grateful for the hospitatility of the Hausdorff Research Institute for Mathematics in Bonn.

\appendix

\section*{Appendix}
\setcounter{section}{1}
\setcounter{equation}{0}

Theorem \ref{th2} for the stretched exponential case relies heavily on the asymptotics for the probability of moderate deviations in \cite{nagaev68}. An English translation of this article can be found in {\em Selected Translations in Mathematics, Statistics and Probability} (1973), Volume 11. The purpose of this appendix is to provide the main results of this difficult to access article, as they apply in our model.\\

\noindent
We use the same notation as in Section 3, introduced in (\ref{notation}), and write the stretched exponential tail of the law $\P$ of the independent integer random variables $\eta_j,\ j=1,2,\ldots,L$ as
\[
\PP\big[\eta_j=n\big]\sim Ae^{-\gamma n^{1-\l}} \qquad \text{as } n\to\infty\ ,
\]
where $A,\gamma$ are positive constants and $0<\l<1$. 
We are interested in the asymptotic behaviour of the probabilities $P_L (N)=\PP\big[ S_L =N\big]$
as $L\to\infty$, where $N=\rho_c L+k$, and $k$ deviates from the typical behaviour $k=O(\sqrt{L})$. This is done in \cite{nagaev68} when $A=\gamma=1$ in a series of five theorems and remarks, where the asymptotics are obtained for five different ranges of $k$-values. In the following, we transcribe these theorems for arbitrary values of $A$ and $\gamma$, which are applied in the proofs of Theorems \ref{th3} and \ref{th4} with $A=A_\lambda$ and $\gamma =b/(1-\lambda)$.\\

\noindent
We will denote by $\theta_L$ a sequence that increases to infinity arbitrarily slowly. Let $r_L(m)$ be a sequence such that
\bea
L\sum_{n=r_L(m)}^{\infty} n^me^{-\gamma n^{1-\l}} \longrightarrow 0 \ \text{ as } L\to\infty\ .
\label{cramer0}
\eea
Note that we may take $r_L(m)=\big(\beta\log L\big)^{\frac{1}{1-\l}}$, for any $\beta>\gamma^{-1}$. We will use $\l^{[t]}(z)$ to denote the first $t$ terms of the Cram\'er series (see e.g. \cite{gnedenko}, Chapter 8), where
\begin{equation}
t=\Big[\frac{1}{\l}\Big]-1\qquad\mbox{and}\qquad \l^{[t]}(z)=\l_0+\l_1z+\cdots+\l_{t-1}z^{t-1}\ .
\label{t}
\end{equation}
In particular, $\l^{[t]}(z)\equiv 0$ when $\lambda> 1/2$. The coefficients $\l_j$ depend on the cumulants of the distribution. Finally, we define
\[
c_\l=(1+\l)(2\l)^{-\frac{\l}{1+\l}}\gamma^{\frac{1}{1+\l}}\ .
\]

\underline{{\em Case $1^o$}}
\begin{equation}
\delta\sqrt{L} <k< (c_\l-\delta)(\sigma^2 L)^\frac{1}{1+\l}\ ,
\label{c1}
\end{equation}
where $\delta>0$ is any sufficiently small fixed number.

\begin{theorem2}
If $k$ is as in (\ref{c1}), then
\[
P_L(N)=\frac{1}{\sigma\sqrt{2\pi L}}\exp\Big\{-\frac{k^2}{2L\sigma^2}+\frac{k^3}{L^2}\l^{[t]}\Big(\frac{k}{L}\Big)\Big\}\big(1+o(1)\big)\ \qquad \text{as }L\to\infty\ .
\]
\label{t1}
\end{theorem2}
\begin{remark}
Under the conditions of Theorem \ref{t1}
\begin{equation}
\frac{\PP\big[S_L=N\big]}{\PP\big[S_L=N;\ \eta_j<r_L(m),\ 1\le j\le L\big]}\longrightarrow 1\ ,
\label{rem1}
\end{equation}
where $m$
can be taken as the smallest positive integer such that $\displaystyle N\left(\frac{k}{N}\right)^{m+1}=o(1)$.
\end{remark}

\underline{{\em Case $2^o$}}
\begin{equation}
(c_\l+\delta)(\sigma^2 L)^\frac{1}{1+\l}<k<\frac{L^{1/{2\l}}}{\theta_L}\ .
\label{c2}
\end{equation}

\begin{theorem2}
If $k$ is as in (\ref{c2}), then
\begin{equation}
P_L(N)=\frac{AL}{\sqrt{1-\frac{\sigma^2\gamma\l(1-\l)L}{k^{1+\l}(1-\alpha)^{1+\l}}}}  \exp\Big\{-\frac{\alpha^2k^2}{2\sigma^2L}+\frac{\alpha^3k^3}{L^2}\l^{[t]}\Big(\frac{\alpha k}{L}\Big) -\gamma(1-\alpha)^{1-\l}k^{1-\l}\Big\}\big(1+o(1)\big)
\label{Pc2}
\end{equation}
as $L\to\infty$, where $\alpha$ is the smallest positive root of the equation
\begin{equation}
\frac{\sigma^2 L}{k^{1+\l}}=\frac{\alpha(1-\alpha)^\l}{\gamma(1-\l)}\big(1-R_\l(\frac{\alpha k}{L})\big)\ .
\label{alef}
\end{equation}
The term $R_\l$ in the preceding equation is given by
\[
R_\l(x)=\frac{\sigma^2}{x}\frac{d}{dx}\big(x^3\l^{[t]}(x)\big)=\sigma^2\sum_{j=0}^{t-1} \l_j(j+3)x^{j+1}\ ,
\]
and in particular $R_\l\equiv 0$ if $\lambda>1/2$.
\label{t2}
\end{theorem2}
\noindent
Note that the fraction of excess mass on the condensate in Theorem \ref{th3} is $a_L =1{-}\alpha$ with $\gamma =b/(1{-}\lambda )$.\\
Relation (\ref{Pc2}) takes a particularly simple form for $\l>1/2$ and $\theta_LL^{\frac{1}{1+\l}}<k<\frac{L^{\frac{1}{2\l}}}{\theta_L}$. In this case
\[
P_L(N)=AL \exp\Big\{-\frac{\alpha^2k^2}{2\sigma^2L}-\gamma(1-\alpha)^{1-\l}k^{1-\l}\Big\}\big(1+o(1)\big).
\]
\begin{remark}
Under the conditions of Theorem \ref{t2}
\begin{equation}
\frac{\PP\big[S_L=N\big]}{L\PP\big[S_L=N;\ \eta_j<r_L(m),\  1\le j \le L-1;\ |\eta_L-(1-\alpha)k|<\theta_L\sqrt{L}\big]}\longrightarrow 1,
\label{rem2}
\end{equation}
where $m$
can be taken as the smallest positive integer such that $\displaystyle Lk^{-\l(m+1)}=o(1)$ as $L\to \infty$.
\end{remark}

\underline{{\em Case $3^o$}}
\begin{equation}
k>L^{\frac{1}{2\l}}\theta_L.
\label{c3}
\end{equation}

\begin{theorem2}
If $k$ is as in (\ref{c3}), then
\[
P_L(N)=AL\exp\big\{-\gamma k^{1-\l}\big\}\big(1+o(1)\big) \qquad \text{as } L\to\infty.
\]
\label{t3}
\end{theorem2}
\begin{remark}
Under the conditions of Theorem \ref{t3}
\[
\frac{\PP\big[S_L=N\big]}{L\PP\big[S_L=N;\ \eta_j<r_L(2),\ 1\le j\le L-1;\ |\eta_L-k|<\theta_L\sqrt{L}\big]}\longrightarrow 1\qquad \text{as }L\to\infty.
\]
\end{remark}

\underline{{\em Case $4^o$}}
\begin{equation}
(c_\l-\delta)(\sigma^2 L)^\frac{1}{1+\l}<k<(c_\l+\delta)(\sigma^2 L)^\frac{1}{1+\l}.
\label{c4}
\end{equation}

\begin{theorem2}
If $k$ is as in (\ref{c4}), then
\begin{align}
P_L(N)=&\frac{AL}{\sqrt{1{-}\frac{\sigma^2\gamma\l(1-\l)L}{k^{1+\l}(1-\alpha)^{1+\l}}}}  \exp\Big\{ {-}\frac{\alpha^2k^2}{2\sigma^2L}{+}\frac{\alpha^3k^3}{L^2}\l^{[t]}\Big(\frac{\alpha k}{L}\Big) {-}\gamma(1{-}\alpha)^{1-\l}k^{1-\l}\Big\}\big(1{+}o(1)\big)\nonumber \\
 &+\frac{1}{\sigma\sqrt{2\pi L}}\exp\Big\{-\frac{k^2}{2L\sigma^2}+\frac{k^3}{L^2}\l^{[t]}\Big(\frac{k}{L}\Big)\Big\}\big(1+o(1)\big) \qquad \text{as }L\to\infty\ ,
\label{Pc4}
\end{align}
where, as before, $\alpha$ is the smallest positive root of equation (\ref{alef}).
\label{t4}
\end{theorem2}
\noindent
Relation (\ref{Pc4}) takes a simple form for $\l>1/2$. In this case
\begin{align}
P_L(N)=&\frac{1}{\sigma\sqrt{2\pi L}}\exp\Big\{-\frac{k^2}{2L\sigma^2}\Big\}\big(1+o(1)\big)\nonumber\\
&+\frac{AL}{\sqrt{1-\frac{\sigma^2\gamma\l(1-\l)L}{k^{1+\l}(1-\alpha)^{1+\l}}}} \exp\Big\{-\frac{\alpha^2k^2}{2\sigma^2L}-\gamma(1-\alpha)^{1-\l}k^{1-\l}\Big\}\big(1+o(1)\big).
\label{simple}
\end{align}
\begin{remark}
This case is intermediate between cases 1 and 2. Under the conditions of Theorem \ref{t4} we can only state that
\[
\frac{\PP\big[S_L=N\big]}{L\PP\big[S_L{=}N;\, \eta_j<r_L(m),\ j\le L{-}1;\, |\eta_L{-}(1{-}\alpha)k|<\theta_L\sqrt{L}\big]+\PP\big[S_L {=}N;\, \eta_j<r_L(m),\forall j\big]}\longrightarrow 1\ ,
\]
where $m=t+2$.
\end{remark}

\underline{{\em Case $5^o$}}
\begin{equation}
\frac{L^{\frac{1}{2\l}}}{\theta_L}<k<L^{\frac{1}{2\l}}\theta_L.
\label{c5}
\end{equation}

\begin{theorem2}
If $k$ is as in (\ref{c3}), then
\[
P_L(N)=AL\exp\Big\{-\gamma k^{1-\l}+\frac{L(1-\l)^2\sigma^2}{2k^{2\l}}\Big\}\big(1+o(1)\big) \qquad \text{as } L\to\infty\ .
\]
\label{t5}
\end{theorem2}
\begin{remark}
In this case the picture is the same as in Remark 3.
\end{remark}

\end{document}